\documentclass[11pt]{article}

\usepackage[utf8]{inputenc}
\usepackage{amsmath}
\usepackage{amsthm}
\usepackage{mathtools}
\usepackage{bm}
\usepackage{algorithm}
\usepackage{algpseudocode}
\usepackage{url}
\usepackage[pdftex]{graphicx}
\usepackage{graphicx}
\usepackage{booktabs}
\usepackage[all]{xy}
\usepackage{mathrsfs}

\usepackage{url}
\usepackage{amsfonts}
\usepackage{fullpage}
\usepackage{natbib} 

\newtheorem{Definition}{Definition}[section]
\newtheorem{Theorem}{Theorem}[section]
\newtheorem{Proposition}{Proposition}[section]
\newtheorem{Lemma}{Lemma}[section]
\newtheorem*{Proof}{Proof}

\newtheorem{Corollary}{Corollary}[section]

\theoremstyle{remark}

\newcommand{\beq}{\begin{equation}}
\newcommand{\beqa}{\begin{eqnarray}}
\newcommand{\eeq}{\end{equation}}
\newcommand{\eeqa}{\end{eqnarray}}
\newcommand{\non}{\nonumber}
\newcommand{\fr}[1]{(\ref{#1})}

\newcommand{\E}{\mathrm{E}}

\newcommand{\N}{\mathrm{N}}
\newcommand{\Z}{\mathrm{Z}}
\newcommand{\sign}{\mathrm{sign}\,}
\newcommand{\sol}{\mathrm{sol}\,}

\newcommand{\cA}{{\cal A}}
\newcommand{\cC}{{\cal C}}

\newcommand{\cI}{{\cal I}}
\newcommand{\cK}{{\cal K}}
\newcommand{\cH}{{\cal H}}
\newcommand{\cL}{{\cal L}}
\newcommand{\cM}{{\cal M}}
\newcommand{\cO}{{\cal O}}

\newcommand{\cS}{{\cal S}}

\newcommand{\dr}{\mathrm{d}}
\newcommand{\wt}[1]{\widetilde{#1}}

\newcommand{\ol}[1]{\overline{#1}}
\newcommand{\ii}{\imath}

\newcommand{\mbbP}{\mathbb{P}}
\newcommand{\mbbR}{\mathbb{R}}

\newcommand{\leadstoup}[1]{\overset{\mathrm{#1}}{\leadsto}}

\begin{document}

\title{Fast symplectic integrator for Nesterov-type acceleration method}

\author{Shin-itiro Goto\\
  Center of Mathematics for Artificial Intelligence and Data Science,\\ 
  Chubu University, \quad 
  1200 Matsumoto-cho, Kasugai, Aichi 487-8501, Japan,\quad
  and\\       
Hideitsu Hino \\
Department of Statistical Modeling, 
The Institute of Statistical Mathematics,\\  
Tachikawa, Tokyo 190-8562, and 
RIKEN AIP, Nihonbashi, Tokyo 103-0027, Japan
}


\maketitle

\begin{abstract}
\end{abstract}
In this paper, explicit stable integrators based on symplectic and contact geometries are proposed for a non-autonomous ordinarily differential equation (ODE) found in improving convergence rate of Nesterov’s accelerated gradient method. Symplectic geometry is known to be suitable for describing Hamiltonian mechanics, and contact geometry is known as an odd-dimensional counterpart of symplectic geometry. Moreover, a procedure, called symplectization, is a known way to construct a symplectic manifold from a contact manifold, yielding Hamiltonian systems from contact ones. It is found in this paper that a previously investigated non-autonomous ODE can be written as a contact Hamiltonian system. Then, by symplectization of a non-autonomous contact Hamiltonian vector field expressing the non-autonomous ODE, novel symplectic integrators are derived. Because the proposed symplectic integrators preserve hidden symplectic and contact structures in the ODE, they should be more stable than the Runge–Kutta method. Numerical experiments demonstrate that, as expected, the second-order symplectic integrator is stable and high convergence rates are achieved.


\section{Introduction}
Optimization plays a central role in solving various engineering problems such as machine learning for data analysis.
To develop effective optimization methods, theoretical proposals that yield fast convergence rates and the design of 
stable numerical schemes for implementations are needed\,\citep{Nocedal2006,Boyd2004,8903465}.

For unconstrained smooth convex problems, several algorithms have been proposed, and Nesterov's accelerated gradient descent (NAG) 
algorithm has been recognized as a major milestone\,\citep{Nesterov1983}.
Subsequently, 
a number of improvements have then been proposed\,\citep{Defazio2019,Odonoghue2015,Fazlyab2008,pmlr-v70-hu17a,NIPS2016_6267,Lessard2016,Hedy2019}. 
One of them derives 
a non-autonomous second-order ordinary differential
equation (ODE) corresponding to the NAG algorithm\,\citep{Su2016}.
Here, an autonomous ODE is synonymous with 
an ODE which does not explicitly depend on an independent variable, where such a variable in this paper is time $t$. 
This second-order ODE enables 
convergence rates and some related quantities to be estimated
if the objectives are sufficiently smooth and the step sizes
in the numerical schemes are sufficiently small.
Determining such an appropriate step size is non-trivial in general. 
That is, there are nontrivial discrepancies
between the continuous-time and discrete-time theories 
of finding an ideal scheme, and several approaches
exist to address this issue. 

For some discrete systems, discrete Lyapunov functions can be
found so that the derivation of a convergence rate does not rely
on the continuous-time limit. Unfortunately, finding
such Lyapunov functions for discrete algorithms
is highly nontrivial\,\citep{bof2018lyapunov,Shi2019,NIPS2019_9508}.
In contrast, links between discrete systems and continuous ones have
been studied in another context.
Because a large number of physical phenomena are modeled
as continuous-time dynamical systems, the construction of
accurate numerical integrators that are suitable 
for continuous-time systems is of great importance.
One sophisticated class of such integrators forms
a class of geometric integrators, and these integrators
preserve the mathematical or geometric structures of the 
continuous-time dynamical system under consideration\,\citep{Hairer2000}.     
In particular {\it{symplectic integrators}} have been applied 
to various Hamiltonian systems including celestial ones because 
symplectic integrators preserve 
the symplectic property in 
Hamiltonian systems, and numerical errors 
do not significantly accumulate 
in long-time simulations\,\citep{Yoshida1990,Yoshida1993}.
Symplectic integrators 
can be employed not only for problems in celestial mechanics,
but also for engineering problems.
In particular, if ODEs that are designed to achieve fast convergence rates
are written as Hamiltonian systems,
then the implementation of symplectic integrators
is expected to realize more stable algorithms than Runge-Kutta integrators
at a fixed order.
To implement symplectic integrators to realize stable algorithms,   
a key mathematical tool should be a procedure for obtaining
Hamiltonian systems from given ODEs.  
Note that there exists 
a systematic procedure yielding a class of {\it{autonomous}} Hamiltonian systems
from a particular class of autonomous ODEs 
by introducing another degree of freedom. 
By contrast, for {\it{non-autonomous}} ODEs, such a procedure is unknown.
This particular class, where autonomous Hamiltonian systems are obtained, is the so-called contact Hamiltonian systems, and the procedure to obtain Hamiltonian systems from contact ones is called symplectization.

\par 

In the recent literature the so-called {\it{contact integrators}}
have been considered.
The theoretical foundation of these integrators is based
on contact geometry, which 
is often called an odd-dimensional analogue of symplectic
geometry\,\citep{Silva2008}.
The significance of the use of contact integrators has
been shown in systems with  
Newtonian mechanics that have 
time-varying non-conserved forces\,\citep{Bravetti2019Celestial}.
In addition, contact integrators can be
applied to accelerated algorithms\,\citep{Bravetti2019arxiv},
and their applicability in optimization should be further explored,
as well as other related 
integrators\,\citep{Vermeeren_2019,frana2020dissipative,Tao-Ohsawa2020}.

\subsection{Summary of this contribution}
A non-autonomous ODE was proposed by\,\citet{Zhang2018}
as a continuous-time limit of NAG so that a fast convergence rate
was obtained. In this paper, using 
the theory of symplectic and contact geometries,
the following are shown:
\begin{itemize}
  \item
It is shown that this ODE is described by a non-autonomous contact 
Hamiltonian system on an extended contact manifold. 
A symplectization of the contact Hamiltonian system 
is then explicitly derived without any approximation. 
Thus, this ODE is shown to belong to both non-autonomous
contact Hamiltonian  and Hamiltonian systems.
\item
A class of symplectic integrators is constructed 
from the obtained non-autonomous 
Hamiltonian system. Because 
the obtained Hamiltonian 
consists of the sum of integrable Hamiltonians,
an explicit scheme is obtained for lower-order integrators. 
Higher-order integrators can in principle be constructed with 
an existing method systematically if necessary. Such 
higher-order integrators can be used if a higher rate is needed. 
The proposed integrators that preserve the Hamiltonian nature 
should be stable. 
This stability was verified numerically.
\end{itemize}
This contribution 
is the first step toward the conceptualization of
the symplectization of contact manifolds as a tool for constructing
high-performance algorithms in optimization problems. 
From a broad perspective, because 
contact geometry is suitable for describing several algorithms
for finding minimizers of objectives\,\citep{Bravetti2019arxiv} and
symplectic geometry is suitable for constructing integrators,
the combination of these two geometries is expected to be 
a fruitful approach 
in optimization problems.

In the following sections the ideas are emphasized,
and the detailed calculations are shown in the 
appendix.  
Related 
notions in differential geometry are also briefly summarized
in the appendix.  
In Section\,\ref{section-geometry},
basic tools for deriving symplectic integrators are summarized,
and these are used to describe Zhang's equation. 
Using these tools, 
symplectic integrators for Zhang's equation are
explicitly derived in Section\,\ref{section-Zhang-equation}.
These theories are numerically verified in Section\,\ref{section-numerics}.
Finally,  Section\,\ref{section-conclusion} summarizes this work briefly. 
\section{Geometric description of non-autonomous ODEs}
\label{section-geometry}
In this section, tools for describing Hamiltonian and contact Hamiltonian systems are briefly summarized, and they are applied to the equation proposed by\,\citet{Zhang2018} and other ones\,\citep{Su2016,Wibisono2016}. 

\subsection{Systems of ODEs as contact and Hamiltonian systems}
Symplectic and contact geometries are well-developed
research areas in mathematics, and they are closely related to each other.
In this subsection a minimum description of these geometries
is summarized to argue  
properties of existing ODEs proposed
as continuous-time optimization methods.   

Symplectic manifold is an even-dimensional manifold together with
the so-called symplectic structure.
In Euclidean setting this is denoted $(\mbbR^{\,2n},\omega)$ with $n\geq1$.
For this manifold, coordinates are $(q,p)$, where
$q=(q^{\,1},\ldots,q^{\,n})\in\mbbR^{\,n}$, and 
$p=(p_{\,1},\ldots,p_{\,n})\in\mbbR^{\,n}$.
The coordinates $q$ express generalized
positions, and $p$ their conjugate momenta.
Symplectic structure $\omega$ is defined as a $2$-form,
and is then also called symplectic $2$-form, 
written in a standard way as
$$
\omega
=\dr\alpha,\qquad
\alpha
=\sum_{a=1}^{n}p_{\,a}\dr q^{\,a}.
$$
where $\alpha$ is a $1$-form called a Liouville $1$-form, $\dr$ is the
so-called exterior derivative. 
Symplectic manifold is suitable for describing a class of energy conservative dynamical systems. Given a function $\ol{H}$ on $\mbbR^{\,2n}$, canonical equations of motion are 
$$
\dot{q}^{\,a}
=\frac{\partial \ol{H}}{\partial q^{\,a}},\qquad
\dot{p}_{\,a}
=-\frac{\partial \ol{H}}{\partial q^{\,a}},\qquad
a=1,\ldots,n.
$$
Here the function $\ol{H}$ is called Hamiltonian that 
is physically identified with energy function, 
$t\in\cI$ denotes time with some $\cI\subseteq\mbbR$, and 
$\dot{\ }$ differentiation with respect to $t$.
One of the roles of $\omega$ is to provide the signs $\pm$ appearing
in the canonical equations. 
It is known that canonical equations of motion are
written as a vector field on $\mbbR^{\,2n}$ in terms of $\omega$. 
Moreover, the energy conservation law, $\dr \ol{H}/\dr t=0$, can be verified.

In the Euclidean setting a contact manifold 
is an odd-dimensional manifold together with the so-called 
{\it{contact structure}}. This is denoted $(\mbbR^{\,2n+1},\ker\lambda)$
with $n\geq1$. 
For this manifold, coordinates are $(q^{\,0},q,\gamma)$, where
$q=(q^{\,1},\ldots,q^{\,n})\in\mbbR^{\,n}$,
$\gamma=(\gamma_{\,1},\ldots,\gamma_{\,n})\in\mbbR^{\,n}$, and $q^{\,0}\in\mbbR$.
In the case where a contact manifold is
applied to describe a class of dissipative 
dynamical systems, $q$ expresses generalized coordinates, $\gamma$ momenta,
and $q^{\,0}$ an action\,\citep{Bravetti2019arxiv}.
Contact structure $\ker\lambda$ 
is defined by the kernel of a contact $1$-form $\lambda$,
$\ker\lambda=\{\,X\in\mbbR^{\,2n+1}|\lambda(X)=0\}$, and thus
contact structure is a $2n$-dimensional subspace of the ambient
$(2n+1)$-dimensional vector space.  
Here a contact $1$-form is written in a standard way as 
$$
\lambda=\dr q^{\,0}-\sum_{a=1}^{n}\gamma_{\,a}\dr q^{\,a}.
$$
Given a function $\ol{K}$ on $\mbbR^{\,2n+1}$, 
equations of motion analogous to canonical equations in the symplectic case are 
\begin{align*}
 \dot{q}^{\,a}
&=\frac{\partial \ol{K}}{\partial \gamma_{\,a}},\quad
\dot{\gamma}_{\,a}
=-\frac{\partial \ol{K}}{\partial q^{\,a}}
-\gamma_{\,a}\frac{\partial \ol{K}}{\partial q^{\,0}},\qquad
a=1,\ldots,n.
\non\\
\dot{q}^{\,0}
&=-\ol{K}+\sum_{a=1}^{n}\gamma_{\,a}\frac{\partial \ol{K}}{\partial \gamma_{\,a}}.
\end{align*} 
Here the function $\ol{K}$ is called contact Hamiltonian.   
Note that the energy conservation law does not hold in general, 
$\dr \ol{K}/\dr t\neq0$.

Relations between symplectic manifolds
and contact manifolds are often argued in autonomous
systems\,\citep{Libermann1987,Schaft2018entropy}.
Symplectization (or sympletification) of a 
$(2n+1)$-dimensional contact manifold is a  
procedure giving $(2n+2)$-dimensional symplectic
manifold. 
By contrast, such a procedure is unknown for 
non-autonomous systems, and    
then the existing symplectization for the autonomous systems
will be extended for 
non-autonomous systems in this paper. This extended symplectization
will be used to construct symplectic integrators
from non-autonomous contact Hamiltonian systems in the later sections.

In the autonomous case, symplectization is summarized as follows.
Given a prescribed contact manifold
$(\mbbR^{\,2n+1},\ker\lambda)$, 
introduce the new variables 
$p_{\,0}\in\mbbR\setminus\{0\}$ and 
\beq
\gamma_{\,a}
=-\,\frac{p_{\,a}}{p_{\,0}},\quad
a=1,\ldots,n.
\label{gamma-p}
\eeq
In addition, introduce a higher dimensional symplectic manifold
$(\mbbR^{\,2n+2},\omega)$, where $\omega$ is given by $\omega=\dr\alpha$. 
This $\alpha$ is called the Liouville $1$-form, and can be constructed from
the contact $1$-form $\lambda$ with Eq.~ \fr{gamma-p}  
as 
$$
p_{\,0}\lambda
=p_{\,0}\dr q^{\,0}+\sum_{a=1}^{\,n}p_{\,a}\dr q^{\,a}
=\alpha,
$$
where $\lambda$ has been treated as
a $1$-form on the symplectic manifold. Similar treatment is   
implicitly applied throughout.    
Then the system with the Hamiltonian
on $\mbbR^{\,2n+2}$
\beq
\ol{H}(q^{\,0},q^{\,1},\ldots,q^{\,n},p_{\,0},p_{\,1},\ldots,p_{\,n})
\eeq
induces the contact Hamiltonian system with $\ol{K}$ on $\mbbR^{\,2n+1}$.
The relation between $\ol{H}$ and $\ol{K}$ is given by 
\beqa
&&\ol{H}(q^{\,0},q^{\,1},\ldots,q^{\,n},p_{\,0},p_{\,1},\ldots,p_{\,n})
\non\\
&&=-\,p_{\,0}\,\ol{K}(q^{\,0},q^{\,1},\ldots,q^{\,n},\gamma_{\,1},
\ldots,\gamma_{\,n}).
\non
\eeqa
This is proven by comparing equations of motion associated with $\ol{H}$ with
those with $\ol{K}$. 

To describe non-autonomous Hamiltonian systems and contact
Hamiltonian systems, introduce additional space $\cI(\subseteq\mbbR)$
whose coordinate
expresses time $t$ for each. 
The resultant phase spaces in Euclidean setting
are $\mbbR^{2n+2}\times\cI$ and $\mbbR^{2n+1}\times\cI$, respectively.
Given a (non-autonomous) Hamiltonian $H$ on $\mbbR^{\,2n+2}\times\cI$, 
the extended Liouville $1$-form 
and extended symplectic $2$-form on $\mbbR^{\,2n+2}\times\cI$ can be
defined as
$$
\alpha^{\,\E}
=\alpha-H\dr t,
$$
and 
\beq
\omega^{\,\E}
=\dr\alpha^{\,\E}.
\label{extended-omega-by-d-alpha}
\eeq
This $\omega^{\,\E}$ is not the standard symplectic form,
because it is defined on an odd-dimensional manifold. 
The canonical 
equations of motion are
\beq
\dot{q}^{\,a}
=\frac{\partial H}{\partial p_{\,a}},\quad 
\dot{p}_{\,a}
=-\,\frac{\partial H}{\partial q^{\,a}},\quad 
\dot{t}=1,
\label{non-autonomous-Hamiltonian-vector-coordinates-Euclidean}
\eeq
where $a=0,\ldots,n$. 
Note that 
$\dr H/\dr t=\partial H/\partial t\neq 0$ in general. 

In addition, given a (non-autonomous) contact Hamiltonian $K$ on
$\mbbR^{\,2n+1}\times\cI$,
the extended contact $1$-form on $\mbbR^{\,2n+1}\times\cI$
are defined by  
$$
\lambda^{\,\E}
=\lambda+K\,\dr t.
$$
The equations of motion are
\begin{subequations}
  \begin{align}
  \dot{q}^{\,a}
&=\frac{\partial K}{\partial \gamma_{\,a}},\quad 
\dot{\gamma}_{\,a}
=-\frac{\partial K}{\partial q^{\,a}}
-\gamma_{\,a}\frac{\partial K}{\partial q^{\,0}},
\\
\dot{q}^{\,0}
&=-K+\sum_{b=1}^{n}\gamma_{\,b}\frac{\partial K}{\partial \gamma_{\,b}},\quad
\dot{t}=1,
\end{align} 
\label{non-autonomous-contact-vector-coordinates-Euclidean}
\end{subequations}
where $a=1,\ldots,n$. 

In this paper $(\mbbR^{\,2n+2}\times\cI,\omega^{\,\E})$ is referred to as
an extended symplectic manifold,
and $(\mbbR^{\,2n+1}\times\cI,\ker\lambda^{\,\E})$ an extended contact manifold.
They are specified after Hamiltonian and contact Hamiltonian are introduced,
respectively.    
A link between non-autonomous Hamiltonian systems and non-autonomous 
contact Hamiltonian systems is as follows.
\begin{Proposition}
\label{fact-non-autonomous-contact-Hamiltonian-lift-Euclidean}
A non-autonomous Hamiltonian system 
associated with $H$ on  $\mbbR^{2n+2}\times\cI$  
induces a non-autonomous contact Hamiltonian system associated
with $K$ on $\mbbR^{\,2n+1}\times\cI$, 
where the relation between $H$ and $K$ is given by
\beqa
&&H(q^{\,0},q^{\,1},\ldots,q^{\,n},p_{\,0},\ldots,p_{\,n},t)
\non\\
&&=-\,p_{\,0}K(q^{\,0},q^{\,1},\ldots,q^{\,n},\gamma_{\,1},\ldots,\gamma_{\,n},t)
\non
\eeqa
where $\gamma_{\,a}=-p_{\,a}/p_{\,0}$ for $p_{\,0}\neq 0$.
Moreover the relation $\alpha^{\,\E}=p_{\,0}\,\lambda^{\,\E}$ hold.
\end{Proposition}
\begin{Proof}
  The claim on equations of motion is verified by
  substituting
the Hamiltonian $H$ 
into Eq.~\fr{non-autonomous-Hamiltonian-vector-coordinates-Euclidean}.
The calculations are as follows. It immediately follows that 
$$
\dot{p}_{\,0}
=-\,\frac{\partial H}{\partial q^{\,0}}
=p_{\,0}\frac{\partial K}{\partial q^{\,0}},\
\dot{q}^{\,0}
=\frac{\partial H}{\partial p_{\,0}}
=-K+\sum_{b=1}^{n}\gamma_{\,b}\frac{\partial K}{\partial \gamma_{\,b}},
$$
and for $a=1,\ldots,n$, 
$$
\dot{p}_{\,a}
=-\,\frac{\partial H}{\partial q^{\,a}}
=p_{\,0}\frac{\partial K}{\partial q^{\,a}},\quad
\dot{q}^{\,a}
=\frac{\partial H}{\partial p_{\,a}}
=\frac{\partial K}{\partial \gamma_{\,a}}.
$$
Combining these calculations, one has 
$$
\dot{\gamma}_{\,a}
=\frac{\dr}{\dr t}\left(-\,\frac{p_{\,a}}{p_{\,0}}\right)
=-\,\gamma_{\,a}\frac{\partial K}{\partial q^{\,0}}
-\frac{\partial K}{\partial q^{\,a}}.
$$
Hence we have derived Eq.~
\fr{non-autonomous-contact-vector-coordinates-Euclidean}.
The claim on $1$-forms is verified as follows. First,
it follows that 
$$
\alpha^{\,\E}-p_{\,0}\,\lambda^{\,\E}
=(\alpha-H\dr t)-p_{\,0}\,(\lambda+K\dr t).
$$
Then substituting $H=-p_{\,0}K$ into the equation above, one has
$$
\alpha^{\,\E}-p_{\,0}\,\lambda^{\,\E}
=\alpha-p_{\,0}\,\lambda.
$$
The right hand side of the equation above vanishes
due to $\alpha=p_{\,0}\lambda$ as in the case of 
the autonomous symplectization. Hence $\alpha^{\,\E}=p_{\,0}\lambda^{\,\E}$. 
\qed
\end{Proof}
The relation between a non-autonomous contact Hamiltonian system and
its symplectization is summarized as follows:  
$$
\xymatrix@C=40pt@R=13pt{
  (\mbbR^{\,2n+2}\times\cI,\,\omega^{\,\E})
  &(q,p,t)\\
  (\mbbR^{\,2n+1}\times\cI,\,\ker\lambda^{\,\E})\ar[u]_{\mbox{Ext. Symplectization}}
  &(q,\gamma,t)\ar@{|->}[u].
}
$$
\subsection{Zhang's equation as contact and Hamiltonian systems}
To design a high performance solver, \citet{Zhang2018} proposed the autonomous ODE based on NAG, 
\beq
\ddot{x}+\Gamma_{\,1}(t;\sigma)\,\dot{x}+\Gamma_{\,0}(t;\sigma)\nabla f (x)
=0,\qquad
\sigma\geq 2, 
\label{zhang2018-ODE} 
\eeq
where $x\in\mbbR^{\,d}$,
$\nabla =(\partial /\partial x^{\,1},\ldots\partial/\partial x^{\,d})$, 
$\dot{x}=\dr x/\dr t\in\mbbR^{\,d}$, 
$t\in\cI:=\mbbR_{>0}$,  
$f$ a prescribed function on $\mbbR^{\,d}$ as an objective,  
$$
\Gamma_{\,0}(t;\sigma)
:=\sigma^{\,2}\,t^{\,\sigma-2},\ \mbox{and}\ 
\Gamma_{\,1}(t;\sigma)
:=\frac{2\sigma+1}{t}.
$$
In what follows, this time-domain $\cI:=\mbbR_{>0}$
is kept fixed so that there is no singularity in $\cI$.   
In this paper Eq.~\fr{zhang2018-ODE} is called Zhang's equation.
\begin{Proposition}
\label{fact-Zhang-is-non-autonomous-contact-Euclidean}
  Zhang's equation \fr{zhang2018-ODE}
  is expressed as a non-autonomous contact Hamiltonian system.
\end{Proposition}
\begin{Proof}
This proof is a generalization of that found in  ~\citet{Bravetti2019arxiv}.
  Identify $n=d$, 
  $x=(q^{\,1},\ldots,q^{\,d})$ so that
  $q=(q^{\,0},x)\in \mbbR^{\,d+1}$, and choose the contact Hamiltonian
  on $\mbbR^{\,2d+1}\times\cI$ as
\beqa
&&  K^{\,\Z}(q^{\,0},q^{\,1},\ldots,q^{\,d},\gamma_{\,1},\ldots,\gamma_{\,d},t)
\non\\
&&  =\frac{1}{2}\sum_{a=1}^{d}\gamma_{\,a}^{\,2}
  +\Gamma_{\,0}(t;\sigma)\,f(q^{\,1},\ldots,q^{\,d})
  +\Gamma_{\,1}(t;\sigma)\,q^{\,0}.
\non
\eeqa
  Then Eq.~\fr{non-autonomous-contact-vector-coordinates-Euclidean} yields
$
\dot{q}^{\,a}
=\gamma_{\,a},\ 
\dot{\gamma}_{\,a}
=-\,\Gamma_{\,0}\frac{\partial f}{\partial q^{\,a}}-\Gamma_{\,1}\gamma_{\,a},
\ a=1,\ldots,d,
$
from which one has  
\beq
  \ddot{q}^{\,a}
+\Gamma_{\,1}\dot{q}^{\,a}  +\Gamma_{\,0}\frac{\partial f}{\partial q^{\,a}}
  =0,
  \qquad a=1,\ldots,d.
\label{Zhang2018-q-coordinates}
  \eeq
The extended contact manifold is $(\mbbR^{\,2d+1}\times\cI,\ker\lambda^{\,\E,\Z})$
  with $\lambda^{\,\E,\Z}=\lambda+K^{\,\Z}\,\dr t$ and
  $\lambda=\dr q^{\,0}-\sum_{a=1}^{d}\gamma_{\,a}\dr q^{\,a}$. 
  \qed
\end{Proof}
Note that the equation $\dot{q}^{\,0}=\cdots$ does not
contribute to the proof of
Proposition\,\ref{fact-Zhang-is-non-autonomous-contact-Euclidean}. 

Then one has the theorem that 
forms the theoretical foundation
for constructing symplectic integrators.
\begin{Theorem}
\label{fact-zhang-equations-are-Hamilton}
Zhang's equation \fr{zhang2018-ODE} is expressed as a non-autonomous
Hamiltonian system.
\end{Theorem}
\begin{Proof}
  By combining
  Propositions\,\ref{fact-non-autonomous-contact-Hamiltonian-lift-Euclidean}
  and \ref{fact-Zhang-is-non-autonomous-contact-Euclidean}, 
  the statement is proven.
\qed
\end{Proof}
Notice that the details of $\Gamma_{\,0}$ and $\Gamma_{\,1}$ are
not necessary for this proof. This leads to the following:
\begin{Corollary}
ODEs including 
the cases where (i) $\Gamma_{\,0}=1, \Gamma_{\,1}=3/t$, \citep{Su2016}, and 
(ii) $\Gamma_{\,0}=\sigma^{\,2}\,t^{\,\sigma-2}, \Gamma_{\,1}=(\sigma+1)/t$,
\citep{Wibisono2016},
can be expressed as non-autonomous Hamiltonian systems.
\end{Corollary}

The explicit form of the Hamiltonian stated in
Theorem\,\ref{fact-zhang-equations-are-Hamilton}
is obtained by applying
Proposition\,\ref{fact-non-autonomous-contact-Hamiltonian-lift-Euclidean}
to Proposition\,\ref{fact-Zhang-is-non-autonomous-contact-Euclidean} as
\beqa
&&H^{\,\Z\Z}(q^{\,0},q^{\,1},\ldots,q^{\,d},p_{\,0},p_{\,1},\ldots,p_{\,d},t)
\non\\
&&=-p_{\,0}\left[
  \frac{1}{2}\sum_{a=1}^{d}\frac{p_{\,a}^{\,2}}{p_{\,0}^{\,2}}
  +\Gamma_{\,0}\,f(q^{\,1},\ldots,q^{\,d})+\Gamma_{\,1}\,q^{\,0}
  \right].\quad\ 
\label{Zhang-Hamiltonian-total}
\eeqa
This $H^{\,\Z\Z}$ induces $(\mbbR^{\,2d+2}\times\cI,\omega^{\,\E,\Z\Z})$ with 
$\omega^{\,\E,\Z\Z}=\dr(p_{\,0}\lambda^{\,\E,\Z})$. 
The explicit forms of equations of motion are obtained from Eq.~
\fr{non-autonomous-Hamiltonian-vector-coordinates-Euclidean}
and Eq.~
 \fr{Zhang-Hamiltonian-total} as
\beqa
&&\dot{q}^{\,0}
=\frac{1}{2}\sum_{b=1}^{d}\frac{p_{\,b}^{\,2}}{p_{\,0}^{\,2}}
-\Gamma_{\,0}f-\Gamma_{\,1}q^{\,0},\;
\dot{p}_{\,0}
=\Gamma_{\,1}\,p_{\,0},
\non\\
&&\dot{q}^{\,a}
=-\frac{p_{\,a}}{p_{\,0}},\;
\dot{p}_{\,a}
=p_{\,0}\Gamma_{\,0}\frac{\partial f}{\partial q^{\,a}},\quad a=1,\ldots,d.
\non
\eeqa
The last three equations above yield Eq.~ \fr{Zhang2018-q-coordinates}.

Features of the system with $H^{\,\Z\Z}$ deriving
Eq.~ \fr{Zhang2018-q-coordinates} are in order.
The equation $\dot{q}^{\,0}=\cdots$ does not
contribute to Zhang's equation, and a solution to the 
equation $\dot{p}_{\,0}=\cdots$ with the explicit form of $\Gamma_{\,1}(t;\sigma)$ is given by
\beq
p_{\,0}(t;\sigma)
=p_{\,0}(1)\,t^{\,2\,\sigma+1},\quad t\in\cI.
\label{Zhang-p0-explicit-solution}
\eeq
Because
the pair of variables $(q^{\,0},p_{\,0})$ in $H^{\,\Z\Z}$ is redundant for
deriving Eq.~\fr{Zhang2018-q-coordinates},
a reduced Hamiltonian system is desired. 
Moreover, the use of the analytical expression 
\fr{Zhang-p0-explicit-solution} is expected to enhance
accuracy of integrators. 

To derive a reduced system that incorporates the analytical expression,
the following  will be applied to $H^{\,\Z\Z}$.  
\begin{Proposition}
\label{fact-reduced-Hamiltonian-general}
  Let $\check{H}$ be a non-autonomous Hamiltonian on
  $\mbbR^{\,2d}\times\cI$, and $\Gamma$ a function of $t\in\cI$. 
  Moreover, let $H_{1}$  be a non-autonomous Hamiltonians defined on 
$\mbbR^{\,2d+2}\times\cI$, 
  and $H_{2}$ a non-autonomous Hamiltonian on 
  $\mbbR^{\,2d}\times\cI$ such that 
\begin{align}
&  H_{1}(q^{\,0},\ldots,q^{\,d},p_{\,0},\ldots,p_{\,d},t)
\non\\
&\quad =-\,p_{\,0}\left[\,\check{H}(q_{\,1},\ldots q_{\,d},p_{\,1},\ldots,p_{\,d},t)
  +\Gamma(t)q^{\,0}\,\right],
\non\\
&  H_{2}(q^{\,1},\ldots,q^{\,d},p_{\,1},\ldots,p_{\,d},t)
\non\\
&\quad =-\,p_{\,0}^{\,\sol}(t)\,\check{H}(q_{\,1},\ldots q_{\,d},p_{\,1},\ldots,p_{\,d},t)
\non
\end{align}
with $p_{\,0}^{\,\sol}(t)$ 
an explicit solution to $\dot{p}_{\,0}=\Gamma(t) p_{\,0}$. 

Then the equations of motion for $q^{\,1},\ldots,q^{\,d}$ and
$p_{\,1},\ldots,p_{\,d}$ obtained by $H_{1}$ 
are also obtained by $H_{2}$.
In addition, extended manifolds $(\mbbR^{\,2d+2},\omega^{\,\E,1})$ 
and $(\mbbR^{\,2d+1},\omega^{\,\E,2})$ are defined for the systems
with $H_{\,1}$ and $H_{\,2}$. 
\end{Proposition} 
\begin{Proof}
  The canonical equations of motion derived from $H_{1}$ are
\begin{align*}
&  \dot{q}^{\,0}
  =-\check{H}-\Gamma(t)\,q^{\,0},\quad
  \dot{p}_{\,0}
  =\Gamma(t)\,p_{\,0},
\\
&  \dot{q}^{\,a}
  =-\,p_{\,0}\frac{\partial \check{H}}{\partial p_{\,a}},\quad
  \dot{p}_{\,a}
  =p_{\,0}\frac{\partial \check{H}}{\partial q^{\,a}},\quad
  a=1,\ldots,d.
\end{align*}
Similarly the equations of motion from $H_{2}$ are
$$
  \dot{q}^{\,a}
  =-\,p_{\,0}^{\,\sol}(t)\frac{\partial \check{H}}{\partial p_{\,a}},\ 
  \dot{p}_{\,a}
  =p_{\,0}^{\,\sol}(t)\frac{\partial \check{H}}{\partial q^{\,a}},\quad
  a=1,\ldots,d.
$$
  Comparing the derived equations from $H_{1}$ with those from $H_{2}$, 
  one completes the first part of the proof. Then,
  the extended symplectic form for the system with $H_{\,2}$
  is introduced by
  $\omega^{\,\E,2}=\dr\alpha^{\,\E,2}$,
  where $\alpha^{\,\E,2}=\sum_{a=1}^{\,d}p_{\,a}\dr q^{\,a}-H_{\,2}\,\dr t$.
  Similarly one can define $\omega^{\,\E,1}$ for 
  the system with $H_{\,1}$.
  \qed
\end{Proof}

With Proposition\,\ref{fact-reduced-Hamiltonian-general}, 
a desired reduced Hamiltonian system for Zhang's equation
is explicitly obtained.
\begin{Proposition}
\label{fact-decomposition-Zhang-Hamiltonian}
Let $H_{\,K}^{\,\Z},H_{\,V}^{\,\Z}$ and $H^{\,\Z}$ be Hamiltonians
on $\mbbR^{\,2d}\times\cI$ 
such that $H^{\,\Z}=H_{\,K}^{\,\Z}+H_{\,V}^{\,\Z}$ with 
\begin{align*}
  H_{\,K}^{\,\Z}(p_{\,1},\ldots,p_{\,d},t)
  =& -\,\frac{1}{2\,p_{\,0}(t;\sigma)}\sum_{a=1}^{d}p_{\,a}^{\,2},\\
  H_{\,V}^{\,\Z}(q^{\,1},\ldots,q^{\,d},t)
  =& -\,p_{\,0}(t;\sigma)\,
\Gamma_{\,0}(t;\sigma)\,f(q^{\,1},\ldots,q^{\,d}).
\end{align*}
Then the Hamiltonian $H^{\,\Z}$ on $\mbbR^{\,2d}\times\cI$   
yields Eq.~\fr{Zhang2018-q-coordinates}, where $p_{\,0}(t;\sigma)$
has been given by Eq.~\fr{Zhang-p0-explicit-solution}.
Moreover the corresponding extended symplectic manifold is
$(\mbbR^{\,2d}\times\cI,\omega^{\,\E,\Z})$, where  
$\omega^{\,\E,\Z}=\dr\alpha^{\,\E,\Z}$ with
 $\alpha^{\,\E,\Z}=\sum_{a=1}^{d}p_{\,a}\dr q^{\,a}-H^{\,\Z}\,\dr t$. 
\end{Proposition}
\begin{Proof}
Applying Proposition\,\ref{fact-reduced-Hamiltonian-general} 
to Eq.~\fr{Zhang-Hamiltonian-total}, one completes the proof.
  \qed
\end{Proof}

The relations among introduced systems are briefly summarized as
$$
\xymatrix@C=70pt@R=20pt{
K^{\,Z}\ar@{~>}[d]_-{\mbox{Prop.\,\ref{fact-Zhang-is-non-autonomous-contact-Euclidean}}}  
&H^{\,\Z\Z}\ar@{~>}[ld]|{\mbox{Thm.\,\ref{fact-zhang-equations-are-Hamilton}}}
\ar@{~>}[d]^-{\mbox{Prop.\,\ref{fact-reduced-Hamiltonian-general}}}
\ar@{~>}[l]_-{\mbox{Prop.\,\ref{fact-non-autonomous-contact-Hamiltonian-lift-Euclidean}}}\\
\mbox{Zhang's Eq.}
&H^{\,\Z}\ar@{~>}[l]^-{\mbox{Prop.\,\ref{fact-decomposition-Zhang-Hamiltonian}}}
} 
$$
where $A\leadstoup{C}B$ indicates that $A$ induces $B$ with C. 

Features of the system with $H^{\,\Z}=H_{\,K}^{\,\Z}+H_{\,V}^{\,\Z}$
deriving Eq.~\fr{Zhang2018-q-coordinates}
are in order. The redundant variable $q^{\,0}$ in $H^{\,\Z\Z}$
does not appear in $H^{\,\Z}$,
and the analytical solution $p_{\,0}(t)$
is incorporated in $H^{\,\Z}$. 
In addition,  
although the whole system $H^{\,\Z}$ could be non-integrable, 
an analytical solution can be derived
for each system solely with $H_{\,K}^{\,\Z}$ and that with $H_{\,V}^{\,\Z}$
(see Sections\,
\ref{section-flow-H_K} and\,\ref{section-flow-H_V}
in the appendix). 
Unlike this case, there is no such a property for the system with $H^{\,\Z\Z}$.  
As will be discussed in the next section,
these explicit solutions are beneficial in constructing integrators. 
Hence $H^{\,\Z}$, rather than $H^{\,\Z\Z}$, is focused in the following sections.

\section{Symplectic integrators for Zhang's equation}
\label{section-Zhang-equation}
In this section, 
the basic properties of symplectic integrators for the system with $H^{\,\Z}$
are shown, and then corresponding integrators are constructed explicitly.  
\subsection{Basic properties of integrators}

A distinctive performance indicator   
of a numerical solver for optimization problems  
is the convergence rate for a given class of objectives.
In \citep{Zhang2018}, the convergence rate for Eq.~\fr{zhang2018-ODE} with 
$\sigma\geq 2$ was derived as
\beq
|\,f(x(t))-f(x^{\,*})\,|
\leq\cO(\,t^{\,-\,\sigma}\,),
\label{convergence-rate-zhang2018-ODE}
\eeq
where $x^{\,*}$ is the minimizer for a convex objective $f$.
In addition a discretization error was taken 
into account
for the $s$-th order Runge-Kutta
integrator with a carefully chosen step size, so that the right hand side of Eq.~
\fr{convergence-rate-zhang2018-ODE} was replaced with
$\cO(N^{-\,\sigma\,s/(s+1)})$, where $N$ is 
the total number of iterations. 
The rate \fr{convergence-rate-zhang2018-ODE}  is realized with
a higher-order integrator, and can also be applied to
the system with $H^{\,\Z}$ on $\mbbR^{\,2d}\times\cI$. 
This yields the following.
\begin{Proposition}
The convergence rate for
the non-autonomous Hamiltonian system in Proposition\,
\ref{fact-decomposition-Zhang-Hamiltonian} 
is given by Eq.~\fr{convergence-rate-zhang2018-ODE}.
Moreover this rate is obtained in the higher  
order limit of integrators. 
\end{Proposition}

Symplectic integrators are numerical integrators for Hamiltonian systems with 
the property that  symplectic forms are preserved\,\citep{Yoshida1990}. These 
integrators are known to be stable, and  error is not significantly
accumulated.  
Although symplectic integrators
are mostly implemented for  
autonomous systems in the literature,
non-autonomous systems can also be implemented.   
Symplectic integrators can in general be divided into $2$ classes, one is 
explicit, and the other one implicit.   
Explicit integrators are often 
preferable to design high performance algorithms, 
because a process that numerically solves solutions to
algebraic equations is not required for explicit ones.   
The condition when explicit symplectic integrators are realized
is known.   
Consider the case that
a Hamiltonian is of the separable form $H(q,p,t)=H_{K}(p,t)+H_{V}(q,t)$ 
with some $H_{K}$ and $H_{V}$. If the system solely with $H_{K}$ and
that with $H_{V}$ 
are analytically solved, then
explicit symplectic integrators can be constructed. 
Note that $H$ itself needs not be integrable.  

The Hamiltonian $H^{\,\Z}$
in Proposition\,\ref{fact-decomposition-Zhang-Hamiltonian} 
satisfies the condition that explicit symplectic integrators are realized. 
Unlike this, the system with $H^{\,\Z\Z}$ does not satisfy this condition, hence explicit symplectic integrators should be constructed
for the system with $H^{\,\Z}$.  

Before closing this subsection, the definition of
non-autonomous symplectic integrator  
in Euclidean setting is given. 
\begin{Definition}
  Let $(\mbbR^{\,2n+2},\omega)$ be a symplectic manifold, 
  $H$ a Hamiltonian on $\mbbR^{\,2n+2}\times\cI$, 
and  $(\mbbR^{\,2n+2}\times\cI,\omega^{\,\E})$
the extended symplectic manifold with $\cI\subseteq \mbbR$ and $\omega^{\,\E}$ being given by Eq.~ \fr{extended-omega-by-d-alpha}. 
In addition, let $z(t)$ be a point on $\mbbR^{\,2n+2}$, and
$z^{\,\E}(t)$ a point on $\mbbR^{\,2n+2}\times\cI$  such that $z^{\,\E}(t)=(z(t),t)=(q(t),p(t),t)\in \mbbR^{\,2n+2}\times\cI$ expresses an exact solution to the system with $H$ at $t\in\cI$, and $(z^{(s)}(t),t)$ an approximate solution at $t$ labeled by some $s=1,2,\ldots$. If the two conditions, 1. an error between $z(t+\tau)$ and $z^{\,(s)}(t+\tau)$ appears in a $\tau^{\,s}$-term in the Taylor expansion for $|\tau|\ll 1$ provided that $z(t)=z^{\,(s)}(t)$, and 2.  $\omega|_{\,t+\tau}=\omega|_{\,t}$, are satisfied, then a discrete time-evolution algorithm is called an $s$-th order (non-autonomous) symplectic integrator.
\end{Definition}

\subsection{Explicit representation of integrators}
The Hamiltonian $H^{\,\Z}$ can be
split into the two pieces as stated 
in Proposition\,\ref{fact-decomposition-Zhang-Hamiltonian}, 
and the solution to each piece  
can analytically be obtained 
(see Sections\,
\ref{section-flow-H_K} and\,\ref{section-flow-H_V}
in the appendix). 
This yields exact relations 
between 2 
points on $\mbbR^{\,2d}\times\cI$ as follows. 
\begin{Lemma}
\label{fact-Hamiltonian-split}
The non-autonomous Hamiltonian system 
solely with $H_{\,K}^{\,\Z}$ and that solely 
by $H_{\,V}^{\,\Z}$ satisfy the relations  
\beq
H_{\,K}^{\,\Z}:
\left\{
\begin{array}{l}
q^{\,a}(t+\tau)\\
\ =q^{\,a}(t)+\frac{p_{\,a}(t)}{2\sigma p_{\,0}(1)}\left[
\frac{1}{(t+\tau)^{\,2\,\sigma}}-\frac{1}{t^{\,2\,\sigma}}
\right],
p_{\,a}(t+\tau)
=p_{\,a}(t),
\end{array}
\right.
\label{flow-H_K}
\eeq
and
\beq
H_{\,V}^{\,\Z}
:\left\{
\begin{array}{l}
q^{\,a}(t+\tau)
=q^{\,a}(t),
\\
p_{\,a}(t+\tau)
=p_{\,a}(t)
+\frac{\sigma p_{\,0}(1)}{3}\left[
(t+\tau)^{\,3\,\sigma}-t^{\,3\,\sigma}
\right]\frac{\partial f}{\partial q^{\,a}}(t),
\end{array}
\right.
\label{flow-H_V}
\eeq
where $a=1,\ldots,d$, and $\tau>0$ is constant.
\end{Lemma}  
\begin{Proof}
To verify the relations for $H_{\,K}^{\,\Z}$, integrate
$$
  \dot{q}^{\,a}
  =\frac{\partial H_{\,K}^{\,\Z}}{\partial p_{\,a}}
  =-\frac{p_{\,a}}{p_{\,0}(t;\sigma)},\quad\mbox{and}\quad
 \dot{p}_{\,a}
 =-\frac{\partial H_{\,K}^{\,\Z}}{\partial q_{\,a}}
 =0,
$$
where $a=1,\ldots,d$, over time from $t$ to $t+\tau$. This yields Eq.~ \fr{flow-H_K}.
Moreover, for 
$H_{\,V}^{\,\Z}$, integrate
\beqa
\dot{q}^{\,a}
&=&\frac{\partial H_{\,V}^{\,\Z}}{\partial p_{\,a}}
=0,\quad\mbox{and}
\non\\
\dot{p}_{\,a}
&=&-\frac{\partial H_{\,V}^{\,\Z}}{\partial q_{\,a}}
=p_{\,0}(t;\sigma)\Gamma_{\,0}(t;\sigma)\frac{\partial f}{\partial q^{\,a}},
\non
\eeqa
where $a=1,\ldots,d$. This yields Eq.~\fr{flow-H_V}. See 
\ref{section-flow-H_K} and\,\ref{section-flow-H_V} 
of the 
appendix 
for full derivations.
\qed
\end{Proof}

From Proposition\,
\ref{fact-decomposition-Zhang-Hamiltonian} and
an existing literature~\citep{Suzuki1993,Bravetti2019Celestial}, 
one has the following: 
\begin{Theorem}
\label{fact-proposed-lower-symplectic-integrators}
The concatenations $\wt{\Phi}_{\,H^{\Z},\tau}^{\,1}$ and 
$\wt{\Phi}_{\,H^{\Z},\tau}^{\,2}$ defined below are 
$1$st and $2$nd order symplectic integrators for Eq.~\fr{Zhang2018-q-coordinates}:
 \beqa
&&\wt{\Phi}_{\,H^{\Z},\tau}^{\,1}
:
\wt{z}(t+\tau)
\non\\
&&=\Phi_{X_{\,t},(\tau/2)}\circ\Phi_{\,H_{V}^{\Z},\tau}\circ\Phi_{\,H_{K}^{\Z},\tau}\,
\circ\Phi_{X_{\,t},(\tau/2)}
\wt{z}(t)
\non\\
&&\wt{\Phi}_{\,H^{\Z},\tau}^{\,2}
:
\wt{z}(t+\tau)
\non\\
&&=\Phi_{\,X_{\,t},\tau/2}
\circ\Phi_{\,H_{K}^{\Z},\tau/2}
\circ\Phi_{\,H_{V}^{\Z},\tau}
\circ\Phi_{\,H_{K}^{\Z},\tau/2}
\circ\Phi_{\,X_{\,t},\tau/2}
\,\wt{z}(t),
\non
\eeqa
where $\Phi_{\,X_{\,t},\tau}$ is the time-shift transform $t\mapsto t+\tau$,  
$\Phi_{\,H_{V}^{\Z},\tau}$ and $\Phi_{\,H_{K}^{\Z},\tau}$ are the transforms $z(t)\mapsto z(t+\tau)$ by Eqs.~ \fr{flow-H_V} and \fr{flow-H_K}, respectively.
\end{Theorem}
\begin{Proof}
These integrators are obtained from Lemma\,\ref{fact-Hamiltonian-split} 
via the splitting method in \citep{Suzuki1993}. 
\qed
\end{Proof}
A higher order integrator requires more gradient evaluations of $f$, and its computational load is high in general. For $\wt{\Phi}_{\,H^{\Z},\tau}^{\,2}$, it requires 1 for each iteration, and is focused below.

An explicit representation of the symplectic integrator for $H^{\,\Z}$ is 
as follows.

For 
the $2$nd order integrator $\wt{\Phi}_{\,H^{\Z},\tau}^{\,2}$,
the following transforms
are concatenated:
\beqa
\Phi_{X_{\,t},(\tau/2)}
\left(
\begin{array}{c}
q^{\,a}\\
p_{\,a}\\
t
\end{array}
\right)
&=&\left(
\begin{array}{c}
q^{\,a}\\
p_{\,a}\\
t+\frac{\tau}{2}
\end{array}
\right),
\non\\
\Phi_{\,H_{K}^{\Z},\tau/2}
\left(
\begin{array}{c}
q^{\,a}\\
p_{\,a}\\
t
\end{array}
\right)
&=&\left(
\begin{array}{c}
q_{\,\star}^{\,a}\\
p_{\,a}\\
t
\end{array}
\right),
\non\\
\Phi_{\,H_{V}^{\Z},\tau}
\left(
\begin{array}{c}
q^{\,a}\\
p_{\,a}\\
t
\end{array}
\right)
&=&\left(
\begin{array}{c}
q^{\,a}\\
p_{\,a}^{\,\ast}\\
t
\end{array}
\right),
\non
\eeqa
where $a=1,\ldots,d$, and
\beqa
q_{\,\star}^{\,a}
&=&q^{\,a}
+\frac{p_{\,a}}{2\,\sigma\,p_{\,0}(1)}
  \left[(t+\tau/2)^{-\,2\sigma}-t^{-\,2\sigma}\right],
\non\\
p_{\,a}^{\,\ast}
&=&p_{\,a}
+\frac{\sigma\,p_{\,0}(1)}{3}
\left[(t+\tau)^{\,3\sigma}-t^{\,3\sigma}\right]\,(\nabla f)_{\,a}.
\non
\eeqa

Note that higher-order symplectic integrators
for this system can systematically be obtained with the
method in~\citep{Suzuki1993,Hatano2005}.

As mentioned in~\citep{Zhang2018,Betancourt2018arxiv},
symplectic integrators are known to be stable. In addition,
as shown in
Proposition\,\ref{fact-decomposition-Zhang-Hamiltonian},
Zhang's equation can be written as  
a non-autonomous Hamiltonian system.  
Hence, by combing these, the proposed 
symplectic integrators are expected to be stable.

\section{Numerical experiments}
\label{section-numerics}
In this section the performance of the proposed $2$nd order symplectic integrator~(SI2) is compared with that of the following existing methods
\begin{enumerate}
\item
the $s$-th order Runge-Kutta~(RK) methods for Zhang's equation, where $s=2,4$~\citep{Griffiths2010}; 
\item
the NAG method, described by
\beqa
x^{\,(k)}
&=&y^{\,(k-1)}-s_{\,\N}\,(\nabla f)(y^{\,(k-1)}),
\non\\
y^{\,(k)}
&=&x^{\,(k)}+\frac{k-1}{k+2}\,(x^{\,(k)}-x^{\,(k-1)}), 
\non
\eeqa
where $s_{\,\N}>0$ is a step size parameter,
$x^{\,(k)}$ and $y^{\,(k)}$ 
denote $x\in\mbbR^{\,d}$ and $y\in\mbbR^{\,d}$ at discrete step $k\geq 0$,
respectively.
\end{enumerate}

We consider two-class classification problems by the regularized logistic regression of the form $f(x)=\frac{1}{|D|}  \sum_{i \in D} \ell_i(x)+\lambda_{\,\mathrm{reg}}\,\|x\|_{\,2}^{\,2}$, where $D$ is the training dataset, $\ell_i(x)$ is a logistic loss function for parameter $x$ given the $i$-th datum in $D$ (pair of observed input and output), $\| x\|_{\,2}^{\,2}:=x_{\,1}^{\,2}+\cdots + x_{\,d}^{\,2}$, and $\lambda_{\,\mathrm{reg}}$ is kept fixed with $10^{-8}$ throughout for simplicity.
To this end, the four popular datasets ({\tt{breastcancer, diabetis, housevotes, sonar}}) are chosen  from UCI machine learning repository, and MNIST dataset where the problem is discriminating even and odd numbers~\footnote{All of the experiments are conducted with MacBook Pro with 2.4Ghz 8-core Intel Core i9 and 64GB RAM. Source code to reproduce the experimental results is available from\\
\url{https://github.com/hideitsu/Contact_SymplecticIntegrator}}.  
Because SI2, RK2 and RK4 are derived from the same ODE, the convergence behavior and computational cost are focused and classification accuracy is not discussed. Accuracies of the classifier obtained by SI2 and NAG were almost the same.


\begin{figure*}[t!]
\centering
\includegraphics[scale=0.23]
{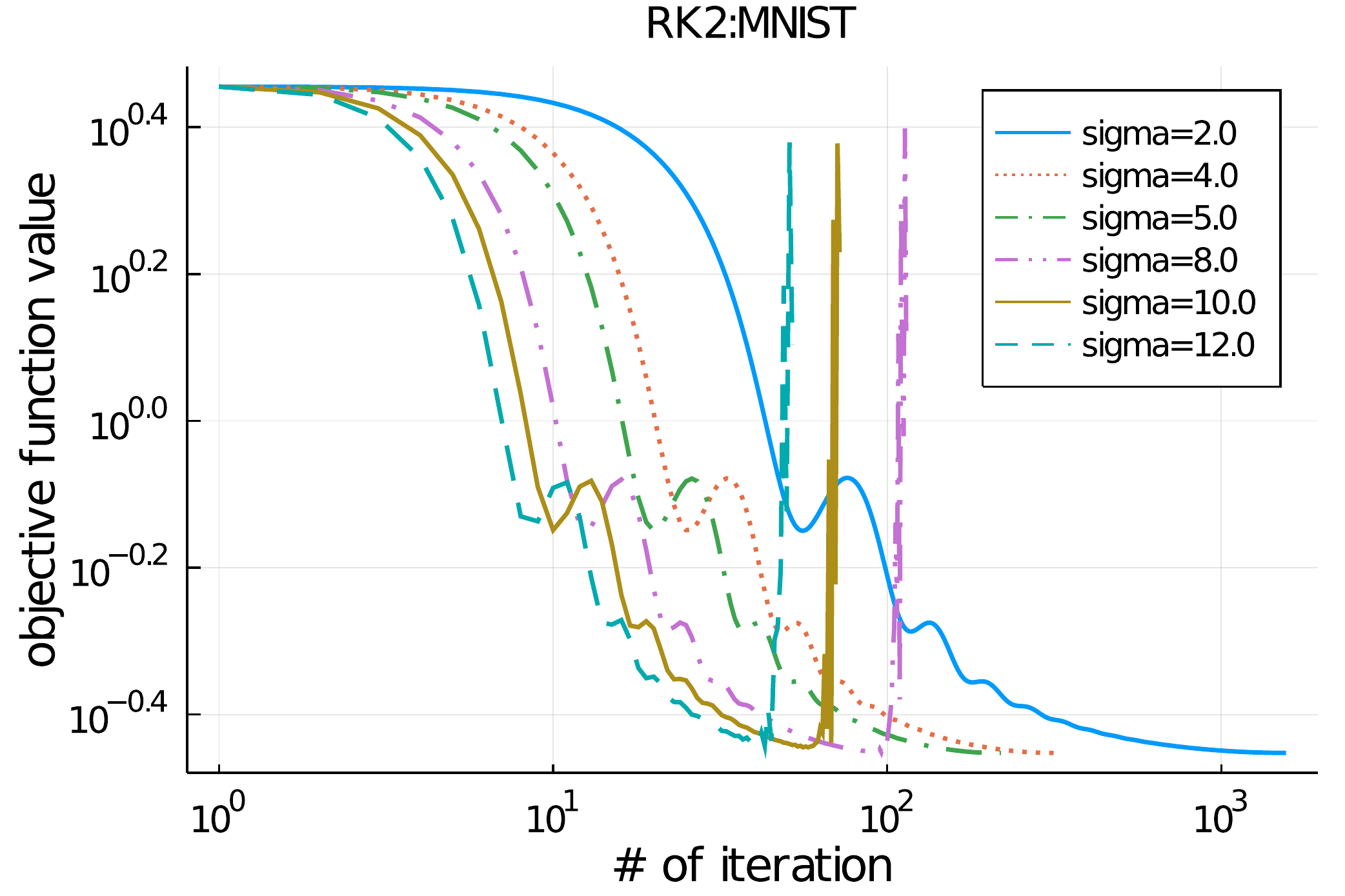}
\includegraphics[scale=0.23]
{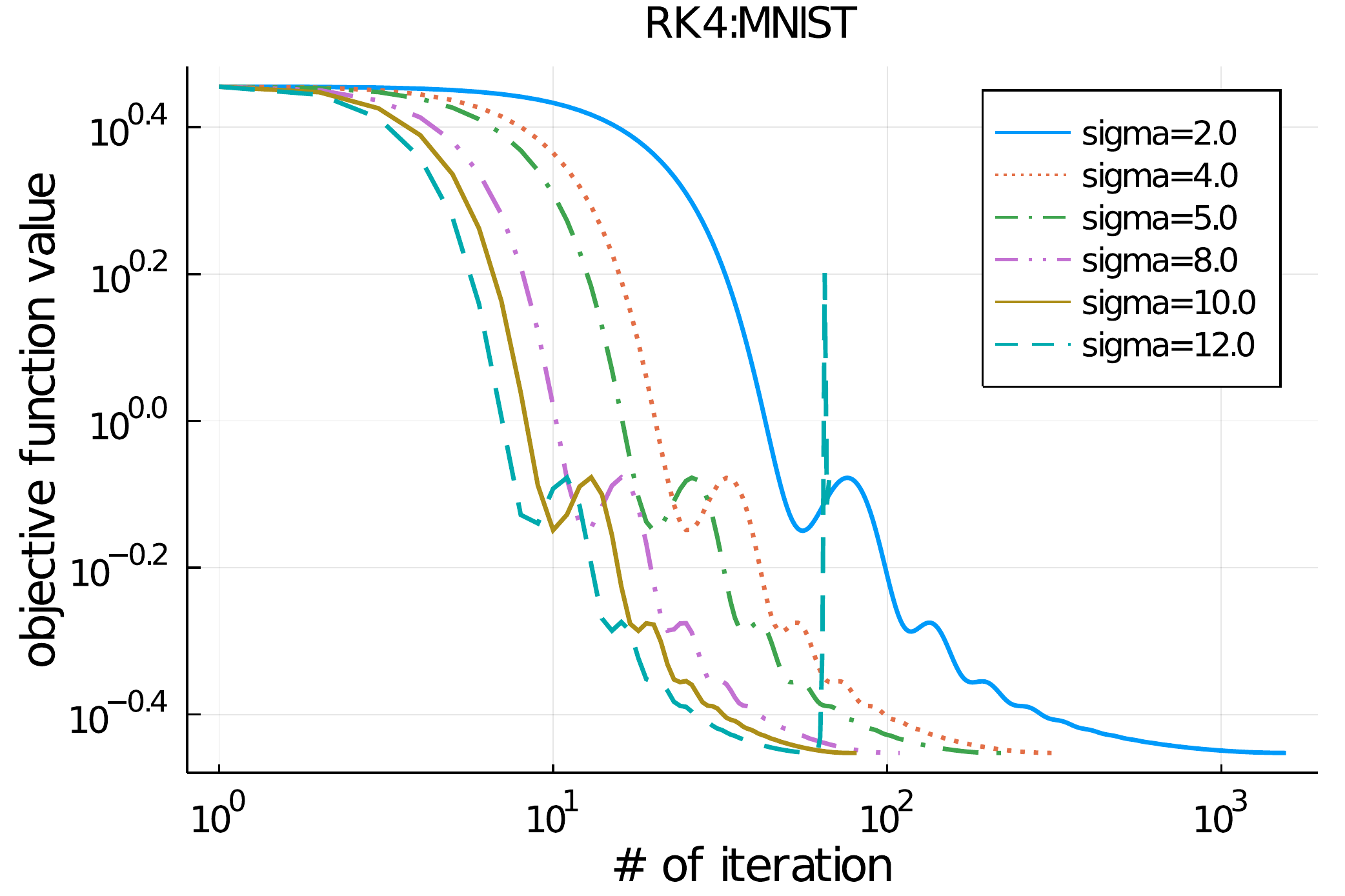}
\includegraphics[scale=0.23]
{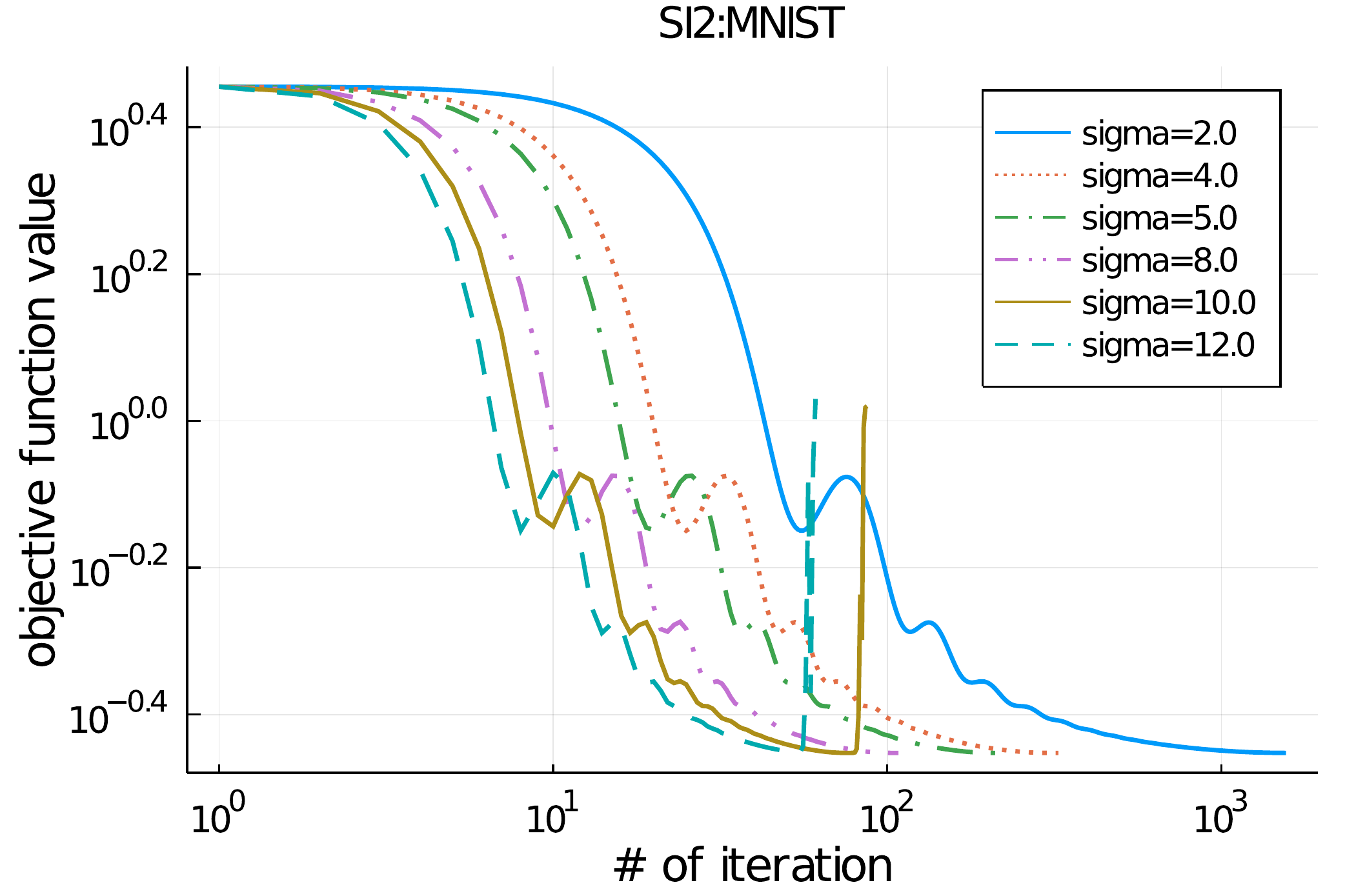}
\caption{Objective values along with the iteration of optimization with varying convergence parameter values of $\sigma$.}
\label{fig:convergence}
\end{figure*}

\begin{table*}[t]
\caption{Computational time (msec) of various methods, averages and SD of 10 trials from different random initialization.}
\label{tab:time}
\vskip 0.15in
\begin{center}
\begin{small}
\begin{sc}
\scalebox{0.8}[0.8]{
\begin{tabular}{lccc|cc}
\toprule
Data set & RK2 & RK4 & SI2 & SI2(bt) & NAG(bt)\\
\midrule
BreastCancer    & 14.51$\pm$ 5.15 & 23.56 $\pm$ 7.89 & 10.85 $\pm$ 5.76 &10.34 $\pm$ 22.658 & 20.41 $\pm$ 19.57\\
Diabetis & 5.96 $\pm$ 4.68 & 7.7 $\pm$ 4.2 & 6.17 $\pm$ 5.75 & 8.24 $\pm$ 27.33 & 14.77 $\pm$ 20.6 \\
HouseVote    & 6.96 $\pm$ 4.67 & 6.9 $\pm$ 0.49 & 5.9 $\pm$ 4.96 &11.49 $\pm$ 26.84 & 9.86 $\pm$ 20.369 \\
Sonar    & 22.9 $\pm$ 6.05 & 38.29 $\pm$ 7.21 & 15.97 $\pm$ 6.1 & 12.34 $\pm$ 25.77 & 37.05$\pm$ 26.27\\
MNIST     & 7430.42 $\pm$ 177.6 & 12886.94 $\pm$ 147.01 & 4717.85 $\pm$ 174.83 &981.44$\pm$ 248.15& 6182.98 $\pm$ 89.43\\
\bottomrule
\end{tabular}
}
\end{sc}
\end{small}
\end{center}
\vskip -0.1in
\end{table*}

\subsection{Comparison to Runge-Kutta methods}
We first compare the performance of  
RK2/4 and SI2  
with different parameter values $\sigma$ in the original ODE~\eqref{zhang2018-ODE}, which controls the convergence speed. Theoretically, by increasing the value of $\sigma$ in the continuous-time limit,  faster convergence rates should be achieved. In reality, partly due to discretization errors, too large $\sigma$ may cause numerical instability. We fixed the learning rate $\tau=0.01$ for all of the three methods, and report the convergence behaviors of objective functions for MNIST dataset in Fig.~\ref{fig:convergence}. Figures for other datasets are shown in the 
appendix. 
From Fig.~\ref{fig:convergence}, it is verified that in general, with larger $\sigma$, the convergence speed is high. However, RK2 often exhibits instability, particularly for MNIST dataset, while SI2 and RK4 are relatively stable.
This instability is partly due to the fact that RK2 is  derived based on the lower order Taylor expansion and it is more deviated from the continuous-time system compared to RK4.

Now we compare the computational speed of the three methods. Table~\ref{tab:time} shows the average of the computational time to reach the stopping criterion (relative difference of the objective function value is less than $10^{-6}$) for RK2, RK4, SI2, and SI2 with backtracking and NAG with backtracking for adjusting step size (see next subsection).
It is seen that the speed with SI2 is slightly faster than that with RK2 or on par for the first four datasets, and faster than that with RK4. For MNIST, which is the largest size dataset among five datasets, the speed with SI2 is significantly faster than those with the  other two methods. This faster computational time of SI2 is due to the fact that one step computation of SI2 requires only one gradient evaluation, while RK2 and RK4 require two and four-times gradient evaluations, respectively.

\subsection{Comparison to NAG method}

We then compare SI2 to NAG. From the results of the previous experiment, we see that for SI2, $\sigma$ less than 8.0 offers stable results. In this section, the convergence rate parameter $\sigma$ is fixed to $6.0$. For both SI2 and NAG, step size remains to be a tuning parameter. We adopt the backtracking method for automatically adjusting the step size in each iteration. In addition, we 
implemented the momentum restarting mechanism to NAG for stabilizing the performance. 

From Table~\ref{tab:time}, column SI2(BT) and NAG(BT), it is seen that computational time of SI2 with backtracking is significantly faster than NAG, particularly for MNIST dataset.

The results of these experiments show that the SI2 is a stable integrator with a high convergence rate.
\section{Conclusion}
\label{section-conclusion}
This paper has shown that the symplectization of a contact manifold
can be used for constructing non-autonomous symplectic integrators 
when ODEs are written as non-autonomous contact Hamiltonian systems.  
In particular, for Zhang's equation~\citep{Zhang2018}, which belongs to a class of continuous-time accelerated gradient methods, explicit non-autonomous symplectic integrators have been constructed. 
Because the proposed symplectic integrators preserve hidden geometric structures in Zhang's equation,  
this should improve performance. 
The resultant $2$nd order integrator
for Zhang's equation was then shown to be
as stable as the $4$-th order Runge-Kutta method,
while the convergence rate was unchanged.
The reason why such explicit integrators can successfully
be obtained is that the split Hamiltonians yield integrable systems. 
However, a profound connection between this integrability
and Zhang's equation is not yet apparent. Because 
Zhang's equation is not only a Hamiltonian system
but also a contact Hamiltonian system, 
the benefits obtained from its contact integrators should also
be explored. 
By addressing these 
questions, we believe that 
faster and more stable algorithms will be realized. In a practical aspect, in experiments on logistic regression using real-world datasets, the convergence speed of the proposed SI2 with larger $\sigma$ has been shown to be faster than that of NAG. It would be important to automatically determine an appropriate value of the convergence rate parameter $\sigma$ for the objective function and the given dataset, which is left to future research. In addition, the estimate of the existing theoretical convergence rate without the backtracking method should be extended for the case where the backtracking method is applied. 

\section*{Acknowledgments}
This work was partially supported by the NEDO Grant Number JPNP18002, JST CREST Grant Number JPMJCR1761, JPMJCR2015. The author S.G was partially supported by JSPS (KAKENHI) Grant No. JP19K03635.

\appendix

\section{Geometric description of non-autonomous ODEs}
\label{section-symplectic-contact-geometries}
To describe the geometric nature of continuous-time dynamical systems
expressed as ODEs, some known facts 
about contact and symplectic geometries are summarized~\citep{Libermann1987,Schaft2018entropy,Bravetti2017ann}.
These descriptions 
are used to reveal the nature of ODEs.
Every geometric object is assumed smooth and real throughout.

\subsection{Autonomous systems}
A symplectic manifold is a pair consisting of an even
dimensional manifold and a closed non-degenerate $2$-form.
Roughly speaking, this manifold is a generalization of a
phase space for an autonomous Hamiltonian system. 
Then contact manifolds are roughly speaking
odd-dimensional counterparts of symplectic manifolds.
Let $n$ be an integer with $n\geq 1$.
A $(2n+1)$-dimensional contact manifold is a pair consisting of 
a $(2n+1)$-dimensional manifold $\cC$ 
and a contact structure $\cA$ defined below.      
First, a $1$-form $\lambda$ on $\cC$ is called 
a contact form if the $(2n+1)$-form 
$\lambda\wedge\dr\lambda\wedge\cdots\wedge\dr\lambda$ does not vanish,
where $\dr$ denotes the exterior derivative, 
and $\wedge$ the exterior or wedge product. 
Then a contact structure on $\cC$ is a 
$2n$-dimensional subspace $\cA\subset T\,\cC$ such that
$\cA=\ker\lambda:=\{X\in T\cC\,|\,\lambda(X)=0 \}$, where $T\cC$ denotes 
the tangent bundle of $\cC$, 
$\lambda(X)$
the (duality) pairing between $\lambda$ and $X$, giving a real number.
As shown bellow, typical symplectic and contact manifolds
are constructed from bundles. 

Let $Q$ be an $(n+1)$-dimensional manifold, $q$ its
coordinates with $q=(q^{\,0},\ldots,q^{\,n})$, and $T^{\,*}Q$
its cotangent bundle. 
A point on $T^{\,*}Q$ is expressed as $(q,p)$ in coordinates 
with $p=(p_{\,0},\ldots,p_{\,n})$. 
Let $\mbbP\,(T^{\,*}Q)$ be the projective cotangent bundle,
which is briefly outlined below.
The manifold $\mbbP\,(T^{\,*}Q)$ is the fiber bundle
whose base space is $Q$, and whose fiber at $q\in Q$ is
the projective space,
$\mbbP\,(T_{\,q}^{\,*}Q):=(T_{\,q}^{\,*}Q\setminus\{0\})/\sim$.
Here this $\sim$ is an equivalence relation, and is given as follows.
If two points $p$ and $p^{\,\prime}$ on $T_{\,q}^{\,*}Q$ are
related by $p^{\,\prime}=\zeta\, p$
with some $\zeta\in\mbbR\setminus\{0\}$, then we write
$p^{\,\prime}\sim p$. 
In a neighborhood where $p_{\,0}\neq 0$, the set
$\gamma=(\gamma_{\,1},\ldots,\gamma_{\,n})$ with 
\beq
\gamma_{\,a}
=-\,\frac{p_{\,a}}{p_{\,0}},\qquad a=1,\ldots,n,
\label{gamma-p-0}
\eeq
can be used as a coordinate set. In other neighborhoods similar
coordinates can be found. 
Then a $(2n+1)$-dimensional manifold $\mbbP(T^{\,*}Q)$ is expressed as
$(q,\gamma)$ in coordinates, 
and $\pi: T^{\,*}Q\to \mbbP\,(T^{\,*}Q)$
denotes a projection, $(q,p)\mapsto(q,\gamma)$.

It is known that $T^{\,*}Q$ induces
a $(2n+2)$-dimensional symplectic manifold $(T^{\,*}Q,\omega)$,
where $\omega=\dr\alpha$ is a symplectic $2$-form
with $\alpha$ being a Liouville $1$-form. 
In addition, $\mbbP\,( T^{\,*}Q)$ induces a $(2n+1)$-dimensional 
contact manifold $(\mbbP\,(T^{\,*}Q),\ker\lambda)$.
According to Darboux's theorem~\citep{arnold1989mathematical},
there exist coordinates $p=(p_{\,0},\ldots,p_{\,n})$ for $T_{\,q}^{\,*}Q$ 
and $\gamma=(\gamma_{\,1},\ldots,\gamma_{\,n})$
for $\mbbP\,(T_{\,q}^{\,*}Q)$ such that
$$
\alpha
=p_{\,0}\dr q^{\,0}+\sum_{a=1}^{n}p_{\,a}\dr q^{\,a},\quad\mbox{and}\quad
\lambda
=\dr q^{\,0}-\sum_{a=1}^{n}\gamma_{\,a}\dr q^{\,a},
$$
in a neighborhood where $p_{\,0}\neq 0$ in $T_{\,q}^{\,*}Q$. 
Then, it follows immediately from $\omega=\dr\alpha$ that 
$$
\omega=\sum_{a=0}^{n}\dr p_{\,a}\wedge\dr q^{\,a}, 
$$
and it is verified that 
$$
  \ker\lambda
  =\mathrm{span}\left(
  \check{\Gamma}_{\,1},\ldots,
  \check{\Gamma}_{\,n},
  \check{\Gamma}^{\,1},\ldots,
  \check{\Gamma}^{\,n}  \right),
  $$
$$
  \check{\Gamma}_{\,a}
  :=\gamma_{\,a}\frac{\partial}{\partial q^{\,0}}+\frac{\partial}{\partial q^{\,a}},\quad
  \check{\Gamma}^{\,a}
  :=\frac{\partial}{\partial \gamma_{\,a}}
  $$
  and
$\alpha=p_{\,0}(\pi^{\,*}\lambda)$, where $\pi^{\,*}$ is the pull-back
induced by $\pi$ 
(pull-backs and push-forwards are explained
in Section\,\ref{section-push-forward-pull-back}). 
A symplectic vector field is a vector field
if it preserves a symplectic $2$-form. Then,  
the Hamiltonian vector field $X_{\,\cH}$ associated with a Hamiltonian
  $\cH$ is a unique vector field satisfying $\ii_{\,X_{\cH}}\omega=-\,\dr \cH$, 
  where $\ii_{\,X}$ denotes the interior product with a vector field 
  $X$ acting on forms.
In contrast, 
contact vector field is a vector field if
it preserves contact structure, $\ker\lambda$. 
Then, the contact Hamiltonian vector field associated with
a function $\cK$ called {\it{contact Hamiltonian}}, denoted by $X_{\,\cK}$,
is a unique vector field satisfying both  
$$
\cL_{\,X_{\cK}}\lambda
=g\,\lambda,\qquad\mbox{and}\qquad
\cK=-\,\ii_{\,X_{\cK}}\lambda,
$$
where $g$ is some function 
and $\cL_{\,X}$
denotes the Lie derivative along a vector field $X$. 
Note that there are several conventions of signs in the literature. 

The following diagram shows
how the introduced manifolds and their vector fields are related:
$$
\xymatrix@R=16pt@C=20pt{
T(T^{\,*}Q)\ar[r]^-{\pi_{\,*}}\ar[d]
&T(\mbbP(T^{\,*}Q))\ar[d]\\
T^{\,*}Q\ar[r]^-{\pi}
&\mbbP(T^{\,*}Q)
}, 
$$
where
$\pi_{\,*}$ is the push-forward induced by $\pi$ sending 
$X_{\,\cH}\in T(T^{\,*}Q)$ to $X_{\,\cK}\in T(\mbbP(T^{\,*}Q))$.
Although the
relation between $X_{\,\cH}$ and $X_{\,\cK}$ via $\pi_{\,*}$ is not used
in this paper, its non-autonomous analogue is used. 

\subsection{Non-autonomous systems}
To describe time-dependent Hamiltonian and
contact Hamiltonian vector fields,
define first $\cS^{\,\E}=T^{\,*}Q\times \cI$
and $\cC^{\,\E}=\mbbP\,(T^{\,*}Q)\times\cI$, where the coordinate
for $\cI\subseteq\mbbR$ expresses time $t$. 
The manifold $\cS^{\,\E}$ is referred to 
as an extended phase space in the context of 
analytical mechanics~\citep{hand_finch_1998}, and $\cS^{\,\E}$
together with a symplectic form is referred to as an extended
symplectic manifold in this paper. 
Then $\cC^{\,\E}$ together with a contact structure is
referred to as an extended contact manifold in this paper. 
Similar to the case of analytical mechanics,
the so-called (time-dependent) Hamiltonian $H$ can be defined 
as a function on $\cS^{\,\E}$. 
Likewise, (time-dependent) contact Hamiltonian
$K$ can be defined on $\cC^{\,\E}$.  
Next, 
given $H$ and $K$, 
$1$-forms on $\cS^{\,\E}$ and $\cC^{\,\E}$ 
are defined as 
$$
\alpha^{\,\E}
=\alpha-H\dr t,
\qquad\mbox{and}\qquad
\lambda^{\,\E}
=\lambda+K\dr t,
$$
where $\alpha^{\,\E}$ is known as the Poincar\'{e}-Cartan form.
It is known that non-autonomous Hamiltonian vector field $X_{\,H}^{\,\E}$
is defined as the one satisfying the condition
$$
\ii_{\,X_{H}^{\E}}\dr\alpha^{\,\E}
=0,
$$
so that $\cL_{\,X_{H}^{\E}}\omega^{\,\E}=0$
with $\omega^{\,\E}:=\dr\alpha^{\,\E}$.
This condition and the choice $X_{\,H}^{\,\E}=X_{\,H}+\partial/\partial t$
with $\dr t(X_{\,H})=0$ give
$$
\ii_{\,X_{H}^{\,\E}}\,\omega
=-\,\dr_{\,T^{\,*}Q}\, H,\quad\mbox{and}\quad
X_{\,H}H=0,
$$
where
$$
\dr_{\,T^{\,*}Q}\, H
:=\dr H-\frac{\partial H}{\partial t}\dr t.
$$ 
In Darboux coordinates, $X_{\,H}^{\,\E}$ is expressed as
\beq
\label{non-autonomous-Hamiltonian-vector-coordinates}
X_{\,H}^{\,\E}
=\sum_{a=0}^{n}\left[
  \dot{q}^{\,a}\frac{\partial}{\partial q^{\,a}}
  +  \dot{p}_{\,a}\frac{\partial}{\partial p_{\,a}}
  \right]+\frac{\partial}{\partial t},
\eeq
where $\dot{q}^{\,a}$ and $\dot{p}_{\,a}$ are given by 
\beq
\label{non-autonomous-Hamiltonian-vector-coordinates-canonical}
\dot{q}^{\,a}
=\frac{\partial H}{\partial p_{\,a}},\ \mbox{and}\ 
\dot{p}_{\,a}
=-\frac{\partial H}{\partial p_{\,a}},
\quad a=0,1,\ldots,n.
\eeq 
A full derivation of Eq.~
\fr{non-autonomous-Hamiltonian-vector-coordinates-canonical}
is given in Section\,\ref{section-non-autonomous-Hamiltonian-vector-coordinates-canonical}.  
By identifying $\dot{}=\dr/\dr t$, 
these equations 
are well-known canonical equations of motion. 
One way to define a vector field extended from an autonomous
contact vector field is as follows~\citep{Bravetti2017ann}.
The non-autonomous contact vector field, denoted by $X_{K}^{\,\E}$,
is defined as the one satisfying both
$$
\cL_{\,X_{K}^{E}}\lambda^{\,\E}
=g^{\,\E}\lambda^{\,\E},\qquad\mbox{and}\qquad
K
=-\ii_{\,X_{K}^{\E}}\lambda,
$$
where $g^{\,\E}$ is some function.
Notice that although the pair
$(\cC^{\,\E},\dr\lambda^{\,\E})$ becomes a symplectic manifold, 
the vector field 
$X_{\,K}^{\,\E}$  does not preserve the symplectic
structure $\dr\lambda^{\,\E}$, in the sense
that $\cL_{\,X_{K}^{\E}}\dr\lambda^{\,\E}\neq 0$.
In Darboux coordinates, $X_{\,K}^{\,\E}$ is expressed as
$$
X_{\,K}^{\,\E}
=\dot{q}^{\,0}\frac{\partial}{\partial q^{\,0}}
+\sum_{a=1}^{n}\left[\,
\dot{q}^{\,a}\frac{\partial}{\partial q^{\,a}}
+\dot{\gamma}_{\,a}\frac{\partial}{\partial \gamma_{\,a}}\,
\right]+\frac{\partial}{\partial t},
$$
where we have chosen the scale factor for $t$ to be unity and 
\begin{subequations}
\begin{align}
\dot{q}^{\,0}
=-\,K
+\sum_{a=1}^{n}\gamma_{\,a}\frac{\partial K}{\partial\gamma_{\,a}},\quad
\dot{q}^{\,a}
=\frac{\partial K}{\partial\gamma_{\,a}},
\\
\dot{\gamma}_{\,a}
=-\frac{\partial K}{\partial q^{\,a}}
-\gamma_{\,a}\frac{\partial K}{\partial q^{\,0}},\quad
a=1,\ldots,n.
\end{align}
\label{non-autonomous-contact-vector-coordinates}
\end{subequations}
A full derivation of Eq.~\fr{non-autonomous-contact-vector-coordinates}
is given in
Section\,\ref{section-non-autonomous-contact-vector-coordinates}.

A way to bridge a contact Hamiltonian vector
field on $\mbbP(T^{\,*}Q)$ and a Hamiltonian vector field on $T^{\,*}Q$ 
is known~\citep{Libermann1987,Schaft2018entropy}.
By extending this existing 
method, the following holds.
\begin{Proposition}
\label{fact-non-autonomous-contact-Hamiltonian-lift-general}
A non-autonomous contact Hamiltonian vector field $X_{\,K}^{\,\E}$
associated with $K$ on  $\mbbP(T^{\,*}Q)\times\cI$ is
lifted to $X_{\,H}^{\,\E}$ associated with an appropriate
$H$ on $T^{\,*}Q\times\cI$, that is,
$\pi_{\,*}^{\,\E}X_{\,H}^{\,\E}=X_{\,K}^{\,\E}$
with $\pi_{\,*}^{\,\E}$ being the push-forward induced by  
$\pi^{\,\E}:T^{\,*}Q\times\cI\to \mbbP(T^{\,*}Q)\times\cI$.
\end{Proposition}
\begin{Proof}
In a neighborhood where $p_{\,0}\neq 0$, choose $H$ to be 
\beqa
&&
H(q^{\,0},q^{\,1},\ldots,q^{\,n},p_{\,0},p_{\,1},\ldots,p_{\,n},t)
\non\\
&&=-\,p_{\,0}\,K
  (q^{\,0},q^{\,1},\ldots,q^{\,n},\gamma_{\,1},\ldots,\gamma_{\,n},t),
\ \gamma_{a}
  =-\frac{p_{\,a}}{p_{\,0}}
\non
\eeqa
  where $p$ is determined by Eq.~\fr{gamma-p-0}.
  Then Eq.~\fr{non-autonomous-Hamiltonian-vector-coordinates}
  and Eq.~\fr{non-autonomous-Hamiltonian-vector-coordinates-canonical}
  yield Eq.~\fr{non-autonomous-contact-vector-coordinates}.
  This is verified as
$$
\dot{p}_{\,0}
=-\frac{\partial H}{\partial q^{\,0}}
=p_{\,0}\frac{\partial K}{\partial q^{\,0}},
\ 
\dot{q}^{\,0}
=\frac{\partial H}{\partial p_{\,0}}
=-\,K+\sum_{a=1}^{n}\gamma_{\,a}\frac{\partial K}{\partial \gamma_{\,a}},
$$
and
$$
\dot{p}_{\,a}
=-\frac{\partial H}{\partial q^{\,a}}
=p_{\,0}\frac{\partial K}{\partial q^{\,a}},\ 
\dot{q}^{\,a}
=\frac{\partial H}{\partial p_{\,a}}
=\frac{\partial K}{\partial \gamma_{\,a}},\ a=1,\ldots,n
$$
from which
$$
\dot{\gamma}_{\,a}
=\frac{\dr}{\dr t}\left(-\frac{p_{\,a}}{p_{\,0}}\right)
=-\,\gamma_{\,a}\frac{\partial K}{\partial q^{\,0}}
-\frac{\partial K}{\partial q^{\,a}},\qquad a=1,\ldots,n.
$$
It can also be proven for other neighborhoods.
\qed
\end{Proof}
In this paper a resultant (non-autonomous) Hamiltonian
vector field obtained by the procedure in
Proposition\,\ref{fact-non-autonomous-contact-Hamiltonian-lift-general}
is called a {\it symplectization} of a (non-autonomous)
contact Hamiltonian vector field. Here, roughly speaking,
symplectization is a procedure for obtaining a symplectic manifold
from a lower-dimensional contact manifold. 
In this paper, extended symplectization is the 
procedure for obtaining an extended symplectic manifold from 
a lower-dimensional extended contact manifold, which
is accomplished by $\omega^{\,\E}=\dr\alpha^{\,\E}$
with $\alpha^{\,\E}=p_{\,0}(\pi^{\,\E\,*}\lambda^{\,\E})$.  
The following diagrams show the relations among the introduced manifolds:
$$
\xymatrix@C=20pt@R=16pt{
T(T^{\,*}Q\times\cI)\ar[r]^-{\pi_{\,*}^{\,\E}}\ar[d]
&T(\mbbP(T^{\,*}Q)\times\cI)\ar[d]
\\
T^{\,*}Q\times \cI\ar[r]^-{\pi^{\,\E}}&\mbbP(T^{\,*}Q)\times\cI
},
$$
which in coordinates, 
$$
\xymatrix@C=20pt@R=13pt{
(q,p,t,\dot{q},\dot{p},\dot{t})\ar@{|->}[r]\ar@{|->}[d]
&(q,\gamma,t,\dot{q},\dot{\gamma},1)\ar@{|->}[d]\\
(q,p,t)\ar@{|->}[r]&(q,\gamma,t)
}.
$$
The symplectization for autonomous systems and extended symplectizations
for non-autonomous systems 
are summarized as  
$$
\xymatrix@C=20pt@R=13pt{
(T^{\,*}Q,\omega)
&(T^{\,*}Q\times\cI,\omega^{\,\E})\\
(\mbbP( T^{\,*}Q),\ker\lambda)\ar[u]^{\mbox{Symplectization}}
&(\mbbP(T^{\,*}Q)\times\cI,\ker\lambda^{\,\E}).\ar[u]^{\mbox{Ext. symplectization}}
}
$$

\section{Brief summary of tools in differential geometry}
The natural mathematical language 
for discussing objects on manifolds 
is in terms of differential forms and their
associated objects, 
because these terms reveal properties that do not depend on
any particular coordinate system. 
This language is suitable for describing
Hamiltonian systems and contact Hamiltonian systems, and 
is hence used in this paper.
In this section some of tools used in this paper are summarized\footnote{ 
We refer the reader to related textbooks, for example, 
(i) Mikio Nakahara.{\it Geometry, Topology and Physics; 2nd ed.} CRC Press 2003,
(ii) Sh\-oshichi Kobayashi and Katsumi Nomizu. {\it Foundations of Differential Geometry; Vol. 1} Interscience Publishers, 1963, 
for further details.}.  
\subsection{Vector fields and differential forms}
At a point $p$ of a differentiable $m$-dimensional manifold $\cM$,
a vector space is denoted $T_{\,p}\cM$.
It follows that
$\dim(T_{\,p}\cM)=m$. Let $x=(x^{\,1},\ldots,x^{\,m})$ be a coordinate set
for $\cM$, then the {\it natural basis} 
$$
\left\{\,
\frac{\partial}{\partial x^{\,1}},\ldots,
\frac{\partial}{\partial x^{\,m}}\,\right\}
$$
spans $T_{\,p}\cM$. At a point $p\in \cM$, 
the sum $X=\sum_{a=1}^{m}X^{\,a}(p)\partial/\partial x^{\,a}$
is an element of
$T_{\,p}\cM$, where $\{X^{\,a}\}_{a=1}^{m}$ is a set of functions.
Points on $\cM$ are often denoted in terms of their coordinates, so that
$X=\sum_{a=1}^{m}X^{\,a}(x)\partial/\partial x^{\,a}$.
The set $T\cM:=\bigcup_{p\in\cM}T_{\,p}\cM$ is called a {\it tangent bundle}, and
it follows that $\dim T\cM=2m$. 

Next, 
the dual of $T_{\,p}\cM$ defined at $p\in\cM$
is denoted $T_{\,p}^{\,*}\cM$. 
By definition, an element
$\alpha\in T_{\,p}^{\,*}\cM$ sends a vector $X\in T_{\,p}\cM$ to a real number
$\mbbR$:  
$$
T_{\,p}^{\,*}\cM : T_{\,p}\cM\to\mbbR,\qquad \alpha: X\mapsto \alpha(X),
$$
where $\alpha(X)\in\mbbR$ is called a {\it pairing}
between $\alpha$ and $X$. Pairing 
is written in some textbooks as 
$$
\langle\,,\rangle :
T_{\,p}^{\,*}\cM\times T_{\,p}\cM\ni(\alpha,X)\mapsto
\langle\alpha,X\rangle:=\alpha(X)\in\mbbR.
$$
By the definition of dual space, 
the sum of two elements $\alpha,\beta\in T_{\,p}\cM$, denoted by
$\alpha+\beta$,   
and scalar multiplication $c\,\alpha$ are equipped on $T_{\,p}\cM$ so that 
$$
(\alpha+\beta)(X)
=\alpha(X)+\beta(X),\ \mbox{and}\ 
(c\,\alpha)(X)
=c[\,\alpha(X)\,], 
$$
where $X\in T_{\,p}\cM,\ c\in\mbbR$.
Then the set $T^{\,*}\cM:=\bigcup_{p\in\cM}T_{\,p}^{\,*}\cM$ is called a
{\it cotangent bundle}. 
An element of $T_{\,p}^{\,*}\cM$ is called a {\it $1$-form}, and
one typical example is $\dr f\in T_{\,p}^{\,*}\cM$
where $f$ is a function on $\cM$. 
This $\dr$ is a map 
that sends a function to a $1$-form so that
$$
(\dr f)(X)
=Xf,\qquad
(Xf)(p)
=X_{\,p}f,\qquad p\in\cM.
$$
where a vector field $X$ can act on a function $f$ 
and $Xf$ is another function. 
This operator $\dr$ will be extended to the one acting on 
wider spaces. A function on a manifold is also called a {\it $0$-form}. 
Let $\Lambda^{\,0}\cM$
be a space of 
$0$-forms on $\cM$, and $\Lambda^{\,1}\cM$ 
a space of $1$-forms.
Then, one can write 
$$
\dr :\Lambda^{\,0}\cM\to \Lambda^{\,1}\cM.
$$
There are various bases for $T_{\,p}^{\,*}\cM$, and one of them is
$$
\left\{\ \dr x^{\,1},\ldots,\dr x^{\,m}\ \right\},
$$
so that the pairing between $\dr x^{\,a}$ and $\partial/\partial x^{\,b}$ is  
$$
\dr x^{\,a}\left(\frac{\partial}{\partial x^{\,b}}\right)
=\delta_{\,b}^{\,a}
:=\left\{
\begin{array}{cc}
  1&a=b\\
  0&a\neq b
\end{array}
\right..
$$
In terms of the basis $\{\dr x^{\,a}\}_{a=1}^{m}$ and
a set of functions $\{X_{\,a}\}_{a=1}^{m}$, 
a $1$-form at $p\in\cM$ can be written as the sum 
$\sum_{a=1}^{m}X_{\,a}(p)\dr x^{\,a}$. 
For a given function $f$, the $1$-form
$\dr f$ is calculated from the definition as 
$$
\dr f
=\sum_{a=1}^{m}\frac{\partial f}{\partial x^{\,a}}\,\dr x^{\,a}.
$$
  This expression is derived below.
Because  $\dr f$ is a $1$-form, it can be written as a sum with the basis
  $\{\dr x^{\,a}\}_{a=1}^{m}$ as 
  $$
  \dr f
  =\sum_{a=1}^{m}g_{\,a}\dr x^{\,a},
  $$
with $\{g_{\,a}\}_{a=1}^{m}$ being a set of functions to be determined below.
Let $X$ be a vector written in terms of  
 basis $\{\partial/\partial x^{\,a}\}_{a=1}^{m}$, that is, 
$$
X=\sum_{a=1}^{m}X^{\,a}\frac{\partial}{\partial x^{\,a}},
$$
where $\{X^{\,a}\}_{a=1}^{m}$ is a set of functions.
Then it immediately follows that 
$$
Xf=\sum_{a=1}^{m}X^{\,a}\frac{\partial f}{\partial x^{\,a}},
$$
and
$$
(\dr f)(X)
=\sum_{b=1}^{m}g_{\,b}\,\dr x^{\,b}\left(
\sum_{a=1}^{m}X^{\,a}\frac{\partial}{\partial x^{\,a}}\right)
=\sum_{a=1}g_{\,a}X^{\,a}.
$$
Comparing these two equations, one has
$$
g_{\,a}
=\frac{\partial f}{\partial x^{\,a}},\qquad a=1,\ldots,m.
$$


Given $\alpha=\sum_{a=1}^{m}f_{\,a}\,\dr x^{\,a}\in T_{\,p}^{\,*}\cM$ and
$X=\sum_{b=1}^{m}g^{\,b}\,\partial/\partial x^{\,b}\in T_{\,p}\cM$, the paring
between them is calculated as 
\beqa
\alpha(X)
&=&\sum_{a=1}^{m}f_{\,a}\,\dr x^{\,a}\left(\,
\sum_{b=1}^{m}g^{\,b}\,\frac{\partial}{\partial x^{\,b}}\,\right)
\non\\
&=&\sum_{a=1}^{m}\sum_{b=1}^{m}f_{\,a}\,g^{\,b}
\left[\dr x^{\,a}\left(\,\frac{\partial}{\partial x^{\,b}}\,\right)\right]
\non\\
&=&\sum_{a=1}^{m}\sum_{b=1}^{m}f_{\,a}\, g^{\,b}\delta_{\,b}^{\,a}
\non\\
&=&\sum_{a=1}^{m}f_{\,a}\, g^{\,a}.
\non
\eeqa

A {\it (differential) $k$-form (field)} $\alpha\in\Lambda^{\,k}\cM$
with $0\leq k\leq m$ defines a map at $p\in\cM$ 
$$
\alpha_{\,p}: \underbrace{T_{\,p}\cM\times\cdots\times T_{\,p}\cM}_{k}
\to \mbbR
$$
equipped with the properties
\beqa
\mbox{(i)}&&
\alpha(X_{\,1},\ldots,fX_{\,r}+f^{\,\prime}X_{\,r}^{\,\prime}\ldots,X_{\,k})
\non\\
&&\ 
=f\,\alpha(X_{\,1},\ldots,X_{\,r},\ldots,X_{\,k})
\non\\
&&\qquad +f^{\,\prime}\,\alpha(X_{\,1},\ldots,X_{\,r}^{\,\prime}\ldots,X_{\,k}),
\non\\
\mbox{(ii)}&&
\alpha(X_{\,1},\ldots,X_{\,r}^{\,\prime},\ldots,X_{\,r},\ldots,X_{\,k})
\non\\
&&=-\, \alpha(X_{\,1},\ldots,X_{\,r},\ldots,X_{\,r}^{\,\prime},\ldots,X_{\,k}),
\non\\
\mbox{(iii)}&&
\alpha(X_{\,1},\ldots,\ldots,X_{\,k})
\non\\
&&\quad\mbox{is a differentiable function on $\cM$.}
\non
\eeqa
Here 
$f,f^{\,\prime}\in\Lambda^{\,0}\cM$,
$X_{\,1},\ldots,X_{\,k}\in T_{\,p}\cM$.
Any $k$-forms with $k>m$ are defined such that $\alpha(X)=0$ for all $X$.
Then the abbreviation $\alpha=0$ can be 
used for $\alpha\in\Lambda^{\,k}\cM$ with $k>\dim\cM$.

The {\it exterior product} or {\it wedge product} $\alpha\wedge\alpha^{\,\prime}$
of the two forms 
$\alpha\in\Lambda^{\,k}\cM$ and $\alpha^{\,\prime}\in\Lambda^{\,l}\cM$ is 
such that
\beqa
&&(\alpha\wedge\alpha^{\,\prime})(X_{\,1},\ldots,X_{\,k},X_{k+1},\ldots,X_{\,k+l})
\non\\
&&=\frac{1}{k!\,l!}\sum_{\sigma}(\sign\,\sigma)\,\alpha(X_{\,\sigma(1)},\ldots,X_{\,\sigma(k)})
\non\\
&&\qquad\cdot\,\alpha^{\,\prime}(X_{\,\sigma(k+1)},\ldots,X_{\,\sigma(k+l)}),
\non
\eeqa
where the numerical factor $k!\,l!$ is replaced with $(k+l)!$ in 
another convention. 
Moreover, for the $\sigma$ 
permutation of $(1,\ldots,k+l)$, $\sign\sigma=1$ for even permutations and
$\sign\sigma=-1$ for odd permutations.
The following hold:
\beqa 
\mbox{(i)}&&
\alpha\wedge\alpha
=0,\qquad
\alpha\in\Lambda^{\,k}\cM,\quad
\mbox{where $k$ is odd}, 
\non\\
\mbox{(ii)}&&
\alpha^{\,\prime}\wedge\alpha
=(-1)^{\,k\,l}\,\alpha\wedge\alpha^{\,\prime},
\non\\ 
&&\qquad
\mbox{where $\alpha\in\Lambda^{\,k}\cM,\alpha^{\,\prime}\in\Lambda^{\,l}\cM$},
\non\\
\mbox{(iii)}&&
(\alpha\wedge\alpha^{\,\prime})(X,X^{\,\prime})
=  \alpha(X)\,\alpha^{\,\prime}(X^{\,\prime})
-  \alpha(X^{\,\prime})\alpha^{\,\prime}(X),
\non\\
&&\qquad\mbox{where $\alpha,\alpha^{\,\prime}\in\Lambda^{\,1}\cM$.}
\non
\eeqa 

The definition of the {\it exterior derivative} is then extended to the operator
$$
\dr : \Lambda^{\,k}\cM\to \Lambda^{\,k+1}\cM,\qquad 0\leq k\leq m,
$$
equipped with the properties
\beqa 
\mbox{(i)}&&
  \dr (\,c\,\alpha+c^{\,\prime}\,\alpha^{\,\prime}\,)
  =  c\,\dr \alpha+c^{\,\prime}\,\dr\alpha^{\,\prime},\qquad c,c^{\,\prime}\in\mbbR,
\non\\
\mbox{(ii)}&&
\dr (\alpha\wedge\alpha^{\,\prime})
=(\dr \alpha)\wedge\alpha^{\,\prime}+(-1)^{\,k}\alpha\wedge(\dr\alpha^{\,\prime}),
\non\\
&&\qquad
\mbox{where $\alpha\in\Lambda^{\,k}\cM,
\alpha^{\,\prime}\in\Lambda^{\,l}\cM$}, 
\non\\
\mbox{(iii)}&&
  \dr^{\,2}\alpha
  =0,
  \non\\
&&\qquad  \mbox{where $\alpha\in\Lambda^{\,k}\cM, \dr^{\,2}\alpha:=\dr(\dr\alpha)$,}
\non\\
\mbox{(iv)}&&
  (\dr f)(X)
  =Xf,\qquad X\in T_{\,p}\cM,\quad f\in\Lambda^{\,0}\cM.
\non
\eeqa 

{\it Interior product} $\ii_{\,X}:\Lambda^{\,k}\cM\to\Lambda^{\,k-1}\cM$ with
$X\in T_{\,p}\cM$ is
$$
(\ii_{\,X}\alpha)(X_{\,1},\ldots,X_{\,k-1})
=\alpha(X,X_{\,1},\ldots,X_{\,k-1}),
$$
and
$$
\ii_{\,X}f
=0,\  f\in\Lambda^{\,0}\cM.
$$
    Then it follows that
\beqa 
\mbox{(i)}&&
  \ii_{\,X}(\,f\alpha+f^{\,\prime}\alpha^{\,\prime})
  = f\, \ii_{\,X}\alpha+f^{\,\prime}\ii_{\,X}\alpha^{\,\prime},
  \non\\
&&\qquad\mbox{where $ f,f^{\,\prime}\in\Lambda^{\,0}\cM,\quad \alpha,
  \alpha^{\,\prime}\in\Lambda^{\,k}\cM$},
\non\\
\mbox{(ii)}&&
  \ii_{\,X}(\,\alpha\wedge\alpha^{\,\prime})
  =(\,  \ii_{\,X}\alpha)\wedge\alpha^{\,\prime}
  +(-1)^{\,k}  \alpha\wedge(\,\ii_{\,X}\alpha^{\,\prime}),
  \non\\
&&\qquad\mbox{where $\alpha\in\Lambda^{\,k}\cM,\ 
  \alpha^{\,\prime}\in\Lambda^{\,l}\cM$},
\non\\ 
\mbox{(iii)}&&
  \ii_{\,X}^{\,2}\alpha
  =0,\non\\
  &&\qquad\mbox{where $\alpha\in\Lambda^{\,k}\cM,\quad
  \ii_{\,X}^{\,2}\alpha:=\ii_{\,X}(\ii_{\,X}\alpha)$.}
\non
\eeqa 
Combining $\ii_{\,X}$ and $\dr$, one has that 
$\ii_{\,X}\dr f=\dr f(X)=Xf$ for a function $f$. 
\subsection{Push-forward and pull-back}
\label{section-push-forward-pull-back}
Let $\cM$ and $\cM^{\,\prime}$ be manifolds
whose dimensions need not be the same, 
and let $\varphi:\cM\to\cM^{\,\prime}$ be 
an invertible map. 
This map induces the map 
$$
\varphi_{\,*}:T_{\,p}\cM\to T_{\,\varphi(p)}\cM^{\,\prime},
\qquad p\in\cM
$$
which is called the {\it push-forward} induced by $\varphi$.
Let $x=(x^{\,1},\ldots,x^{\,m})$ be a set of coordinates for $\cM$ with
$\dim\cM=m$, and
$y=(y^{\,1},\ldots,y^{\,n})$ with $\dim\cM^{\,\prime}=n$.
For $X\in T_{\,p}\cM$ given by 
$$
X=\sum_{a=1}^{m}f_{\,a}(x)\frac{\partial}{\partial x^{\,a}},
$$
its push-forward $\varphi_{\,*}$ is calculated to be
$$
\varphi_{\,*}X
=\sum_{a=1}^{m}\sum_{b=1}^{n}
f_{\,a}(x(y))\frac{\partial y^{\,b}}{\partial x^{\,a}}(y)
\frac{\partial}{\partial y^{\,b}},\quad \in T_{\,\varphi(p)}\cM^{\,\prime}.
$$

The map $\varphi:\cM\to\cM^{\,\prime}$  induces another map
$$
\varphi^{\,*} :
T_{\,\varphi(p)}^{\,*}\cM^{\,\prime}\to
T_{\,p}^{\,*}\cM,\qquad p\in\cM,
$$
which is called the {\it pull-back} induced by $\varphi$.
This is also defined
for $k$-forms such that
\beqa
&&(\varphi^{\,*}\alpha)(X_{\,1},\ldots,X_{\,k})
=\alpha(\varphi_{\,*}X_{\,1},\ldots,\varphi_{\,*}X_{\,k}),
\non\\
&&
\qquad X_{\,1},\ldots,X_{\,k}\in T_{\,p}\cM, \quad
\alpha\in \Lambda^{\,k}\cM^{\,\prime}.
\non
\eeqa
More precisely,
\beqa
&&(\varphi^{\,*}\alpha)_{\,p}(X_{\,1},\ldots,X_{\,k})
=\alpha_{\varphi(p)}(\varphi_{\,*}X_{\,1},\ldots,\varphi_{\,*}X_{\,k}),
\non\\
&&
\qquad 
X_{\,1},\ldots,X_{\,k}\in T_{\,p}\cM, \quad
\alpha\in \Lambda^{\,k}\cM^{\,\prime}.
\non
\eeqa
For functions, the pull-back is defined as
$$
\varphi^{\,*}f
=f\circ\varphi,\qquad
f\in\Lambda^{\,0}\cM.
$$
It follows that
\beqa 
\mbox{(i)}&&
\varphi^{\,*}(\alpha\wedge\alpha^{\,\prime})
=(\varphi^{\,*}\alpha)\wedge(\varphi^{\,*}\alpha^{\,\prime}),
\non\\
&&\qquad \alpha\in\Lambda^{\,k}\cM,\quad
\alpha^{\,\prime}\in\Lambda^{\,l}\cM,
\non\\
\mbox{(ii)}&&
  \varphi^{\,*}(\,\dr\alpha\,)
  =\dr  (\,\varphi^{\,*}\alpha\,),\qquad
  \alpha\in\Lambda^{\,k}\cM.
\non
\eeqa 
\subsection{One-parameter group of transforms and the Lie derivative }
The Lie derivative is a tool for evaluating changes in various objects
in differential geometry, and it is used in
Section\,\ref{section-symplectic-contact-geometries}. 
This is briefly explained here.

Suppose that a diffeomorphism $\varphi_{\,t}:\cM\to\cM$
is given for each $t\in\mbbR$. If the map
$$
\varphi : \mbbR\times\cM\to\cM,\quad
\varphi(t,x)
=\varphi_{\,t}(x),\quad
t\in\mbbR,\quad x\in\cM,
$$
is differentiable and satisfies
$$
\varphi_{\,t}\circ\varphi_{\,t^{\,\prime}}
=\varphi_{\,t+t^{\,\prime}},
$$
then $\{\varphi_{\,t};t\in\mbbR\}$ is called a
{\it $1$-parameter  group of transforms}. 
If a domain for $t$ is a bounded domain of $\mbbR$,
$\{\varphi_{\,t}\}$ is called {\it $1$-parameter group of local transforms}. 

Given a $\varphi_{\,t}$, there is a vector field $X$ satisfying 
$$
X_{\,p}f
=\left.\frac{\dr}{\dr t}\right|_{t=0}f(\,\varphi_{\,t}\,p\,),\qquad
p\in\cM,
$$
for any function $f$. This vector field $X$ is called
an {\it infinitesimal transform}.
In contrast, 
given a vector field $X$, there is a
$1$-parameter group of local transforms $\{\Phi_{\,X,t}\}$ such that
$X$ is the infinitesimal transform. 
This induced $\{\Phi_{\,X,t}\}$ is called the {\it $1$-parameter group of (local)
transforms induced by $X$}. The word ``local'' is often omitted.

Let $X$ be a vector field, and $\{\Phi_{\,X,t}\}$
the $1$-parameter group of (local) transform.
Given $\alpha\in\Lambda^{\,k}\cM$, the $k$-form 
$$
\cL_{\,X}\alpha
:=\lim_{t\to 0}\frac{1}{t}\left(\Phi_{\,X,t}^{\,*}\,\alpha-\alpha\right),
$$
is called the {\it Lie derivative }of $\alpha$ along $X$.
Although Lie derivative can be defined for other objects, including 
vector fields, 
such details are not discussed here.
For forms, it follows that
\beqa 
\mbox{(i)}&&
  \cL_{X}(\alpha\wedge\alpha^{\,\prime})
  =  (\cL_{X}\alpha)\wedge\alpha^{\,\prime}
  +  \alpha\wedge(\cL_{X}\alpha^{\,\prime}),
\non\\  
&&  \alpha\in\Lambda^{\,k}\cM,\quad  \alpha^{\,\prime}\in\Lambda^{\,l}\cM,
\non\\
\mbox{(ii)}&&
  \cL_{\,X}\alpha
  =(\ii_{\,X}\dr+\dr\ii_{\,X})\,\alpha,\qquad
  \alpha\in\Lambda^{\,k}\cM,
\non\\
&&  \mbox{(this is know as the Cartan formula)}, 
\non\\
\mbox{(iii)}&&
  \dr\cL_{\,X}\alpha
  =\cL_{\,X}\dr\alpha,\qquad
  \alpha\in\Lambda^{\,k}\cM, 
\non\\
\mbox{(iv)}&&
  \cL_{\,X}f
  =Xf,\qquad
  f\in\Lambda^{\,0}\cM.
\non
\eeqa 

\section{Detailed derivations of equations} 
\subsection{Equation \fr{non-autonomous-Hamiltonian-vector-coordinates-canonical}}
\label{section-non-autonomous-Hamiltonian-vector-coordinates-canonical}
Suppose that a Hamiltonian $H$ is given
on the manifold $\cS^{\,\E}=T^{\,*}Q\times\cI$ with $\dim Q=n+1$.
Then the Poincar\'e-Cartan form
$\alpha^{\,\E}=\alpha-H\,\dr t$ is obtained, which is written in coordinates as 
$$
\alpha^{\,\E}
=\sum_{a=0}^{n}p_{\,a}\dr q^{\,a}-H\dr t,
$$
so that
$$
\omega^{\,\E}
=\dr\alpha^{\,\E}
=\omega-\dr H\wedge\dr t,
\quad\mbox{where}\quad
\omega
=\sum_{a=0}^{n}\dr p_{\,a}\wedge\dr q^{\,a}.
$$
Let $\dot{q},\dot{p},\dot{t}$ be some functions on $\cS^{\,\E}$. 
Then a non-autonomous Hamiltonian vector field with
$$
X_{\,H}^{\,\E\,\prime}
=\sum_{a=0}^{n}\left[\dot{q}_{\,a}\frac{\partial}{\partial q^{\,a}}
+\dot{p}^{\,a}\frac{\partial}{\partial p_{\,a}}
  \right]+\dot{t}\frac{\partial}{\partial t} 
$$
reduces to Eq.~\fr{non-autonomous-Hamiltonian-vector-coordinates}
when the condition $\dot{t}=1$ is imposed. 

A vector field $X_{H}^{\E}\in T\cS^{\,\E}$ satisfying
$$
\ii_{\,X_{H}^{\E}}\omega^{\,\E}
=0 
$$
expresses non-autonomous Hamiltonian equations of motion, that is, 
Eq.~\fr{non-autonomous-Hamiltonian-vector-coordinates-canonical}.
This statement, and the decomposed form 
$$
\ii_{\,X_{H}^{\,\E}}\,\omega
=-\,\dr_{\,T^{\,*}Q}\, H,\quad\mbox{and}\quad
X_{\,H}H=0,
$$
where
$$
\dr_{\,T^{\,*}Q}\, H
:=\dr H-\frac{\partial H}{\partial t}\dr t,
$$
are verified below.

First the decomposed form is obtained. Substituting
\beqa
\ii_{\,X_{H}^{\E}}\omega^{\,\E}
&=&\ii_{\,X_{H}^{\E}}\omega-\ii_{\,X_{H}^{\E}}(\dr H\wedge \dr t)
\non\\
&=&\ii_{\,X_{H}^{\E}}\omega
-(\ii_{\,X_{H}^{\E}}\dr H)\,\dr t
+(\ii_{\,X_{H}^{\E}}\dr t)\,\dr H
\non\\
&=&
\ii_{\,X_{H}^{\E}}\omega
-(X_{\,H}^{\,\E} H)\,\dr t
+(\,X_{\,H}^{\,\E} t)\,\dr H
\non\\
&=&
\ii_{\,X_{H}^{\E}}\omega
-\left[(X_{\,H} H)+\frac{\partial H}{\partial t}\right]\,\dr t
+\dr H
\non\\
&=&
\ii_{\,X_{H}^{\E}}\omega
-(X_{\,H} H)\,\dr t
+\dr_{\,T^{\,*}Q}\, H,
\non
\eeqa
into the condition $\ii_{\,X_{H}^{\E}}\omega^{\,\E}=0$,  
and noticing that the $1$-form 
$$
\ii_{\,X_{H}^{\E}}\omega+\dr_{\,T^{\,*}Q}\, H
$$
does not contain $\dr t$, 
one has the decomposed form.

Second, Eq.~\fr{non-autonomous-Hamiltonian-vector-coordinates-canonical}
is derived below. Because 
\beqa
\ii_{X_{H}^{\E}}\omega
&=&
\sum_{a=0}^{n}\left[(\ii_{X_{H}^{\E}}\dr p_{\,a})\,\dr q^{\,a}
-(\ii_{X_{H}^{\E}}\dr q^{\,a})\,\dr p_{\,a}
\right]
\non\\
&=&\sum_{a=0}^{n}\left[(X_{H}^{\E} p_{\,a})\dr q^{\,a}
-(X_{H}^{\E} q^{\,a})\,\dr p_{\,a}
\right]
\non\\
&=&\sum_{a=0}^{n}\left[\dot{p}_{\,a}\,\dr q^{\,a}
-\dot{q}^{\,a}\,\dr p_{\,a}
\right],
\non\\
-\dr_{\,T^{\,*}Q}\, H
&=&-\dr H+\frac{\partial H}{\partial t}\dr t
\non\\
&=&-\sum_{a=0}^{n}\left[
  \frac{\partial H}{\partial q^{\,a}}\dr q^{\,a}
+\frac{\partial H}{\partial p_{\,a}}\dr p_{\,a}
  \right],
\non
\eeqa
and $\dr q^{\,a},\dr p_{\,a},\dr t$ form a basis on $T^{\,*}S^{\,\E}$, one has
the desired equations
$$
\dr q^{\,a}
:\ 
\dot{p}_{\,a}
=-\,\frac{\partial H}{\partial q^{\,a}},
\, 
\dr p_{\,a}
:\ 
\dot{q}^{\,a}
=\frac{\partial H}{\partial p_{\,a}},\ 
a=1,\ldots,n.
$$
\subsection{Equation \fr{non-autonomous-contact-vector-coordinates}}
\label{section-non-autonomous-contact-vector-coordinates}
Suppose that a contact Hamiltonian $K$ is given on the manifold
$\cC^{\,\E}=\mbbP(T^{\,*}Q)\times\cI$ with $\dim Q=n+1$. Then
$$
\lambda^{\,\E}
=\lambda+K\dr t 
$$
is obtained, where
$\lambda$ is a contact $1$-form on $\cC=\mbbP(T^{\,*}Q)$. In a neighborhood
it is expressed as
$$
\lambda
=\dr q^{\,0}-\sum_{a=1}^{n}\gamma_{\,a}\dr q^{\,a}.
$$

The non-autonomous contact vector field, denoted by $X_{K}^{\,\E}$, 
is defined as the one satisfying both
$$
\cL_{\,X_{K}^{E}}\lambda^{\,\E}
=g^{\,\E}\lambda^{\,\E} \qquad\mbox{and}\qquad
K
=-\ii_{\,X_{K}^{\E}}\lambda,
$$
where $g^{\,\E}$ is some function.
In Darboux coordinates, $X_{\,K}^{\,\E}$ is expressed as
$$
X_{\,K}^{\,\E}
=\dot{q}^{\,0}\frac{\partial}{\partial q^{\,0}}
+\sum_{a=1}^{n}\left[\,
\dot{q}^{\,a}\frac{\partial}{\partial q^{\,a}}
+\dot{\gamma}_{\,a}\frac{\partial}{\partial \gamma_{\,a}}\,
\right]+\frac{\partial}{\partial t},
$$
and $\dot{q}^{\,0}$, $\dot{q}^{\,a}$, and $\dot{\gamma}_{\,a}$ obey  
Eq.~\fr{non-autonomous-contact-vector-coordinates}.

This statement is verified below. 
First, the $2$nd condition is equivalent to
$$
\lambda^{\,\E}(X_{\,K}^{\,\E})
=0,
$$
because 
$$
0=\ii_{\,X_{K}^{\E}}\lambda+K
=\ii_{\,X_{K}^{\E}}(\lambda+K\dr t)
=\ii_{\,X_{K}^{\E}}\lambda^{\,\E}
=\lambda^{\,\E}(X_{\,K}^{\,\E}).
$$
From this $2$nd condition and the Cartan formula, one has that 
$$
\cL_{X_{K}^{\E}}\lambda^{\,\E}
=(\ii_{X_{K}^{\E}}\dr+\dr\ii_{X_{K}^{\E}})\lambda^{\,\E}
=\ii_{X_{K}^{\E}}\dr \lambda^{\,\E},
$$
from which  the $1$st condition can be written as 
$$
\ii_{X_{K}^{\E}}\dr \lambda^{\,\E}
=g^{\,\E}\lambda^{\,\E}.
$$
The right hand side of the equation above is
$$
g^{\,\E}\lambda^{\,\E}
=g^{\,\E}\left(
\dr q^{\,0}-\sum_{a=1}^{n}\gamma_{\,a}\dr q^{\,a}+K\,\dr t
\right).
$$
The left hand side reduces from   
\beqa
\dr\lambda^{\,\E}
&=&\dr\left(\dr q^{\,0}-\sum_{a=1}^{n}\gamma_{\,a}\dr q^{\,a}+K\,\dr t\right)
\non\\
&=&\dr K\wedge \dr t-\sum_{a=1}^{n}\dr\gamma_{\,a}\wedge\dr q^{\,a} 
\non
\eeqa
to  
\beqa
&&\ii_{X_{K}^{\E}}\dr\lambda^{\,\E}
\non\\
&&=\ii_{X_{K}^{\E}}
\left(\dr K\wedge \dr t\right)
-\ii_{X_{K}^{\E}}\left(\sum_{a=1}^{n}\dr\gamma_{\,a}\wedge\dr q^{\,a}\right)
\non\\
&&=
(\ii_{X_{K}^{\E}}\dr K)\, \dr t-(\ii_{X_{K}^{\E}}\dr t)\,\dr K
-\sum_{a=1}^{n}\left[
(\ii_{X_{K}^{\E}}\dr\gamma_{\,a})\,\dr q^{\,a}
- (\ii_{X_{K}^{\E}}\dr q^{\,a})\,\dr\gamma_{\,a}
\right]
\non\\
&&=
(X_{K}^{\E}K)\, \dr t-\,\dr K
-\sum_{a=1}^{n}\left[
(X_{\,K}^{\,\E}\gamma_{\,a})\,\dr q^{\,a}
- (X_{\,K}^{\,\E} q^{\,a})\,\dr\gamma_{\,a}
\right]
\non\\
&&=
(X_{K}^{\E}K)\, \dr t-\,\dr K
-\sum_{a=1}^{n}\left[
\dot{\gamma}_{\,a}\,\dr q^{\,a}
- \dot{q}^{\,a}\,\dr\gamma_{\,a}
\right].
\non
\eeqa
To reduce the equation further, substituting 
\beqa
&&(X_{K}^{\E}K)\, \dr t-\,\dr K
\non\\
&&=(X_{K}^{\E}K)\, \dr t-\bigg[
  \frac{\partial K}{\partial t}\dr t
  +\frac{\partial K}{\partial q^{\,0}}\dr q^{\,0}
+\sum_{a=1}^{n}\left(
\frac{\partial K}{\partial q^{\,a}}\dr q^{\,a}
+  \frac{\partial K}{\partial \gamma_{\,a}}\dr \gamma_{\,a}\right)
\bigg]
\non\\
&&=\left[(X_{K}^{\E}K)-\frac{\partial K}{\partial t}\right]\dr t
-\frac{\partial K}{\partial q^{\,0}}\dr q^{\,0}
-\sum_{a=1}^{n}\left(
  \frac{\partial K}{\partial q^{\,a}}\dr q^{\,a}
+  \frac{\partial K}{\partial \gamma_{\,a}}\dr \gamma_{\,a}\right)
\non 
\eeqa
into the equation above, one has
\beqa
&&\ii_{X_{K}^{\E}}\dr\lambda^{\,\E}
=\left[(X_{K}^{\E}K)-\frac{\partial K}{\partial t}\right]\dr t
-\frac{\partial K}{\partial q^{\,0}}\dr q^{\,0}
\non\\
&&\qquad 
-\sum_{a=1}^{n}\left[
 \left(\dot{\gamma}_{\,a}+\frac{\partial K}{\partial q^{\,a}}\right)\dr q^{\,a}
  +\left(-\dot{q}^{\,a}+  \frac{\partial K}{\partial \gamma_{\,a}}\right)\dr \gamma_{\,a}\right].
\non
\eeqa
Because $\dr q^{\,0},\dr q^{\,a},\dr\gamma_{\,a},\dr t$
form a basis for $T^{\,*}\cC^{\,\E}$, the $1$st condition yields
\beqa
\dr t&:&
(X_{K}^{\E}K)-\frac{\partial K}{\partial t}
=g^{\,\E}K, 
\non\\
\dr q^{\,0}&:&
\qquad -\frac{\partial K}{\partial q^{\,0}}
=g^{\,\E},
\non\\
\dr q^{\,a}&:&
\ \ \dot{\gamma}_{\,a}+\frac{\partial K}{\partial q^{\,a}}
=\gamma_{\,a}g^{\,\E},
\non\\
\dr \gamma_{\,a}&:&
-\dot{q}^{\,a}+\frac{\partial K}{\partial \gamma_{\,a}}
=0.
\non
\eeqa
From these equations, one has
$$
\dot{q}^{\,a}
=\frac{\partial K}{\partial \gamma_{\,a}},\ 
\dot{\gamma}_{\,a}
=-\frac{\partial K}{\partial q^{\,a}}
-\gamma_{\,a}\frac{\partial K}{\partial q^{\,0}},\  a=1,\ldots,n.
$$
Observe that 
$$
\frac{\dr}{\dr t}K
=X_{\,K}^{\,\E}K
=\frac{\partial K}{\partial t}
-\frac{\partial K}{\partial q^{\,0}}K.
$$

The $2$nd condition is equivalent to
\beqa
0&=&\lambda^{\,\E}(X_{\,K}^{\,\E})
\non\\
&=&\left(\dr q^{\,0}-\sum_{b=1}^{n}\gamma_{\,b}\,\dr q^{\,b}+K\,\dr t\right)
(X_{\,K}^{\,\E})
\non\\
&=&X_{\,K}^{\,\E} q^{\,0}-\sum_{b=1}^{n}\gamma_{\,b}\,X_{\,K}^{\,\E} q^{\,b}+K\,X_{\,K}^{\,\E} t
\non\\
&=&
\dot{q}^{\,0}
-\sum_{a=1}^{n}\gamma_{\,a}\,\dot{q}^{\,a}
+K.
\non
\eeqa
Then one has
$$
\dot{q}^{\,0}
=\sum_{a=1}^{n}\gamma_{\,a}\,\frac{\partial K}{\partial\gamma_{\,a}}-K.
$$
\section{Derivations of equations in the main text}
\subsection{Equation (11) 
in the main text}
\label{section-flow-H_K}
Given the Hamiltonian
$$
H_{\,K}^{\,\Z}(q,p,t)
=-\frac{1}{2p_{\,0}(t;\sigma)}\sum_{b=1}^{d}p_{\,b}^{\,2}, 
$$
with 
$$
p_{\,0}(t;\sigma)
=p_{\,0}(1)\,t^{\,2\,\sigma+1}, 
$$
the canonical equations of motion
Eq.~\fr{non-autonomous-Hamiltonian-vector-coordinates-canonical}
are
\beqa
\dot{q}^{\,a}
&=&\frac{\partial H_{\,K}^{\,\Z}}{\partial p_{\,a}}
=-\frac{p_{\,a}}{p_{\,0}(t;\sigma)},
\non\\
\dot{p}_{\,a}
&=&-\frac{\partial H_{\,K}^{\,\Z}}{\partial q_{\,a}}
=0,
\non
\eeqa
for $a=1,\ldots,d$.
The solution to this set of equations is then obtained by integration.
With
\beqa
-\int_{t}^{\,t+\tau}\frac{\dr t^{\,\prime}}{p_{\,0}(t^{\,\prime};\sigma)}
&=&\frac{-1}{p_{\,0}(1)}
\int_{t}^{\,t+\tau}\frac{\dr t^{\,\prime}}{(t^{\,\prime})^{\,2\sigma+1}}
\non\\
&=&\frac{1}{2\,\sigma\,p_{\,0}(1)}\left[\frac{1}{(t+\tau)^{\,2\sigma}}-\frac{1}{t^{\,2\sigma}}\right],
\non
\eeqa
one has
\beqa
&&q^{\,a}(t+\tau)
=q^{\,a}(t)+\frac{1}{2\,\sigma\,p_{\,0}(1)}\left[\frac{1}{(t+\tau)^{\,2\sigma}}-\frac{1}{t^{\,2\sigma}}\right],
\non\\
&&p_{\,a}(t+\tau)
=p_{\,a}(t),
\non
\eeqa
for $a=1,\ldots,d$.
For the sake of completeness, the solution to the system with $H_{\,K}^{\,\Z}$
is obtained as
$$
q^{\,a}(t)
=q^{\,a}(1)+\frac{1}{2\,\sigma\,p_{\,0}(1)}\left[\frac{1}{t^{\,2\sigma}}-1\right], 
p_{\,a}(t)
=p_{\,a}(1),
$$
for $a=1,\ldots,d$.
\subsection{Equation (12)  
in the main text}
\label{section-flow-H_V}
Given the Hamiltonian
$$
H_{\,V}^{\,\Z}(q,p,t)
= -\,p_{\,0}(t;\sigma)\,\Gamma_{\,0}(t;\sigma)\,f(q),
$$
with 
$$
p_{\,0}(t;\sigma)
=p_{\,0}(1)\,t^{\,2\,\sigma+1},\quad
\Gamma_{\,0}(t;\sigma)
=\sigma^{\,2}\,t^{\,\sigma-2}, 
$$
the canonical equations of motion
Eq.~\fr{non-autonomous-Hamiltonian-vector-coordinates-canonical}
are
\beqa
\dot{q}^{\,a}
&=&\frac{\partial H_{\,K}^{\,\Z}}{\partial p_{\,a}}
=0,\non\\
\dot{p}_{\,a}
&=&-\frac{\partial H_{\,K}^{\,\Z}}{\partial q_{\,a}}
=p_{\,0}(t;\sigma)\,\Gamma_{\,0}(t;\sigma)
\frac{\partial f}{\partial q^{\,a}},
\non
\eeqa
for $a=1,\ldots,d$. 
Then the solution to this set of equations is obtained by integration.
With
\beqa
&&\int_{t}^{\,t+\tau}p_{\,0}(t^{\,\prime};\sigma)\,\Gamma_{\,0}(t^{\,\prime};\sigma)
\,\dr t^{\,\prime}
\non\\
&&=\sigma^{\,2}\,p_{\,0}(1)\int_{t}^{\,t+\tau}(t^{\,\prime})^{\,3\sigma-1}
\,\dr t^{\,\prime}
\non\\
&&=\frac{\sigma^{\,2}\,p_{\,0}(1)}{3\sigma}
\left[(t+\tau)^{\,3\sigma}-t^{\,3\sigma}\right],
\non
\eeqa
one has
\beqa
&&q^{\,a}(t+\tau)
=q^{\,a}(t),
\non\\
&&p_{\,a}(t+\tau)
=p_{\,a}(t)
+\frac{\sigma p_{\,0}(1)}{3}
\left[(t+\tau)^{\,3\sigma}-t^{\,3\sigma}\right]\frac{\partial f}{\partial q^{\,a}},
\non
\eeqa
for $a=1,\ldots,d$.  
For the sake of completeness, the solution to the system with $H_{\,K}^{\,\Z}$
is obtained as
\beqa
&&q^{\,a}(t)
=q^{\,a}(1),
\non\\
&&p_{\,a}(t)
=p_{\,a}(1)
+\frac{\sigma p_{\,0}(1)}{3}
\left[t^{\,3\sigma}-1\right]\frac{\partial f}{\partial q^{\,a}},
\non
\eeqa
for $a=1,\ldots,d$.

\section{Numerical evaluation in more detail}
This section is intended to provide more information about the performance of 
the 2nd symplectic integrator (SI2) that has been discussed 
in Section 4,  
that is, the performance of 
the SI2 is compared with those of 
the existing methods RK2, RK4, and NAG in detail.

To this end, recall that one of the significant characteristics in numerical simulations for a method is how many times demanding or heavy calculations is needed.  
For a given dynamical system $\dot{x}=F(x,t)$ with some differentiable function $F:\mbbR^{\,d}\times\cI\to\mbbR^{\,d}$, suppose that we numerically solve for $x(t)\in\mbbR^{\,d}$ with some integrator. One significant computational load is to evaluate $F$ numerically, and the number of  evaluations depends on the algorithm. 
For SI2, it is the total iteration number. For  RK2, it is twice the total iteration number. For RK4, it is 4 times the total iteration number. This means that the computation cost of SI2 is the lowest among these ODE based algorithms. 

In this section, as a realistic problem with convex objective function, we consider two-class classification by regularized logistic regression. 

\subsection{Objective function}
We adopt a standard logistic loss function with regularization to prevent complete separation:
\begin{align}
\notag
    f(w) =& - \frac{1}{|D|} 
    \sum_{i \in D}
    \left\{
    y_{i} \log ( h(x_i;w)) + (1-y_{i}) \log (1-h(x_i;w))
    \right\}\\
    &+ 
    \lambda_{\,\mathrm{reg}} \| w\|_{2}^{2},
\end{align}
where $D$ is the training dataset, and the regularization parameter $\lambda_{\,\mathrm{reg}}$ is fixed to $10^{-8}$ and $h(x;w) = \frac{1}{1+e^{-w^{\top}x}}$. 
Our aim is to develop a novel optimization algorithm based on geometric notion and not to develop a good classifier, hence the regularization parameter is commonly used for all of the algorithms and not tuned. We also note that in the main text of the paper, $x$ is used for parameter while in the above objective function of logistic regression, we follow convention that $w$ is the parameter of the model and $(x_i, y_i)$ are the observation where $x_i \in \mathbb{R}^{d}$ and $y_{i} \in \{+1,-1\}$. 

\subsection{Datasets}
We use four popular datasets for classification from UCI machine learning repository, and the famous MNIST dataset. Profile of these datasets are summarized in table~\ref{tab:data}.

\begin{table*}[t]
\caption{Profile of datasets.}
\label{tab:data}
\vskip 0.15in
\begin{center}
\begin{small}
\begin{sc}
\begin{tabular}{lccccc}
\toprule
Data set & BreastCancer & Diabetis & HouseVote & Sonar & MNIST\\
\midrule
Dimension    & 9& 8 & 16 & 60 & 784\\
Sample size & 615 & 353 & 209 & 187 & 60000\\
\bottomrule
\end{tabular}
\end{sc}
\end{small}
\end{center}
\vskip -0.1in
\end{table*}

\subsection{Comparison to Runge-Kutta methods}
\begin{figure*}[t!]
\centering
\includegraphics[scale=0.23]
{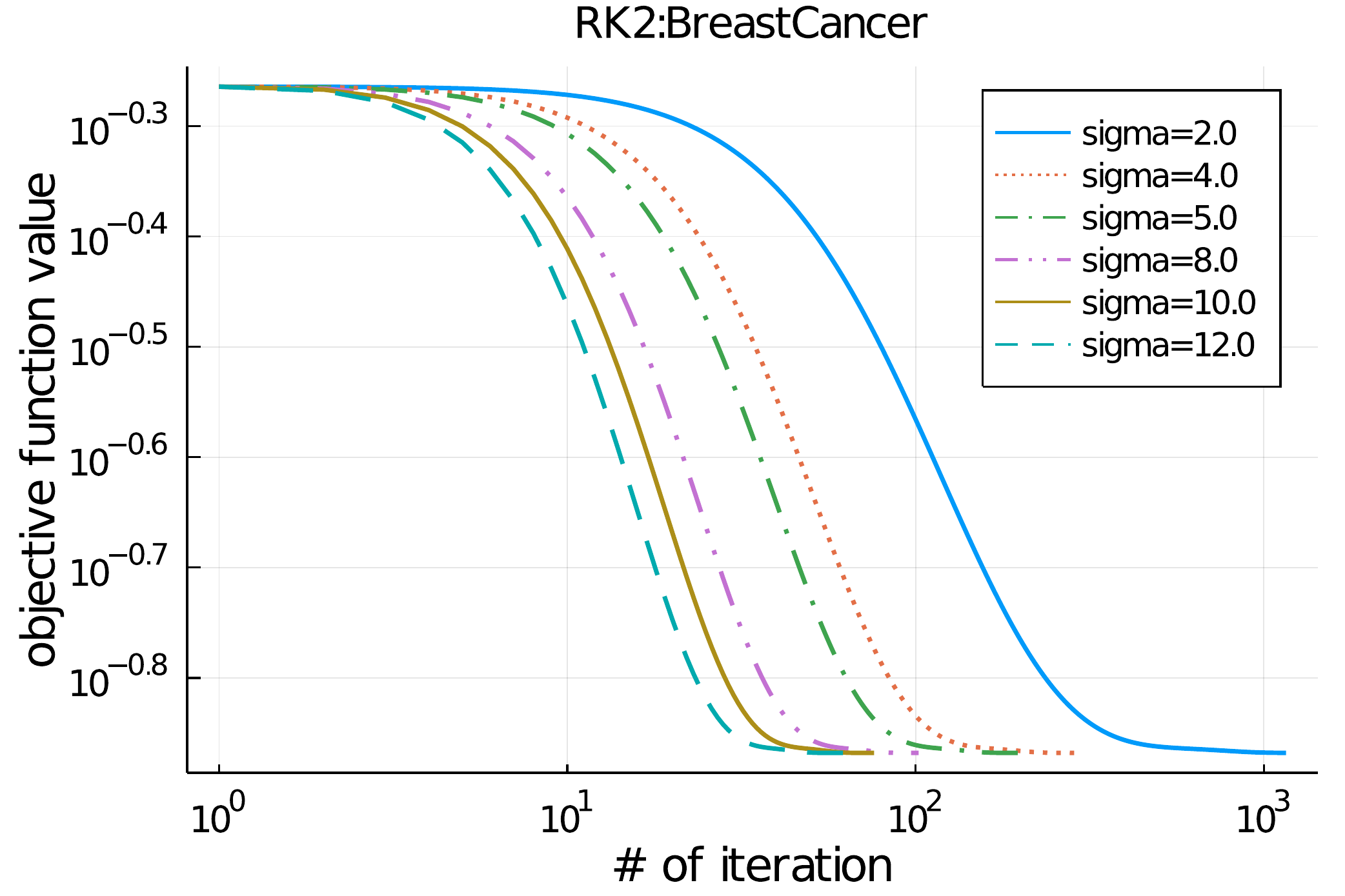}
\includegraphics[scale=0.23]
{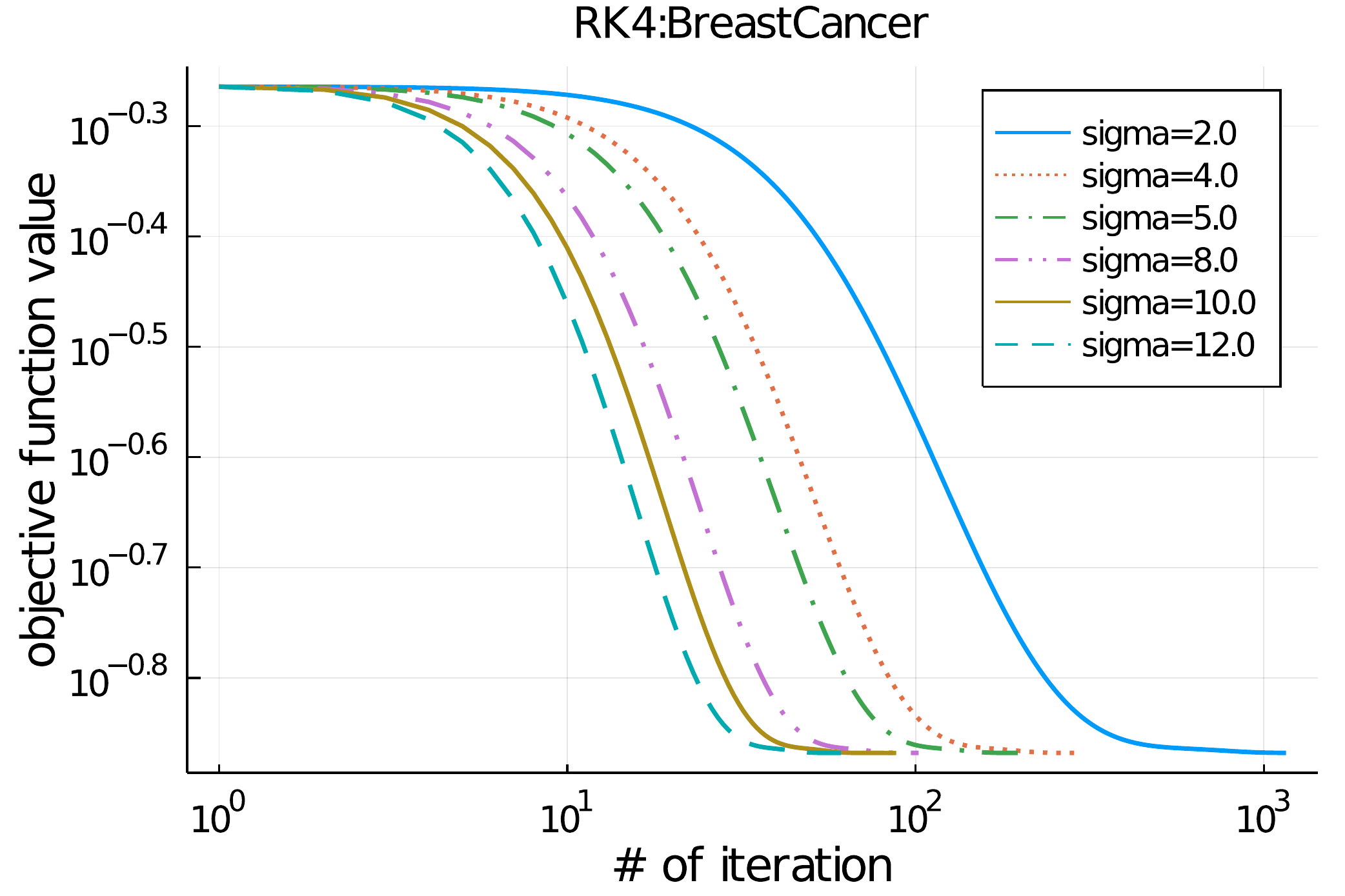}
\includegraphics[scale=0.23]
{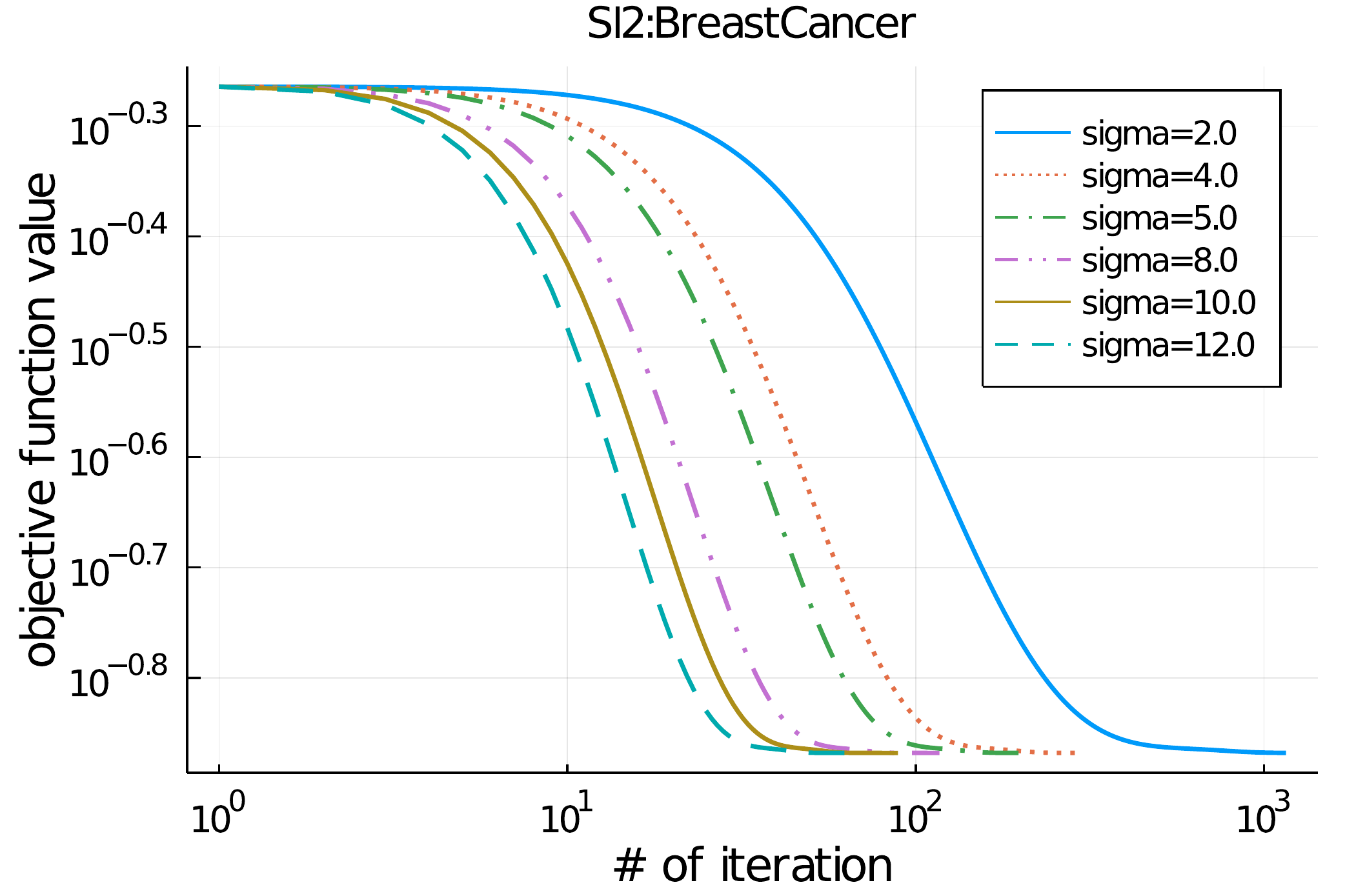}\\
\includegraphics[scale=0.23]
{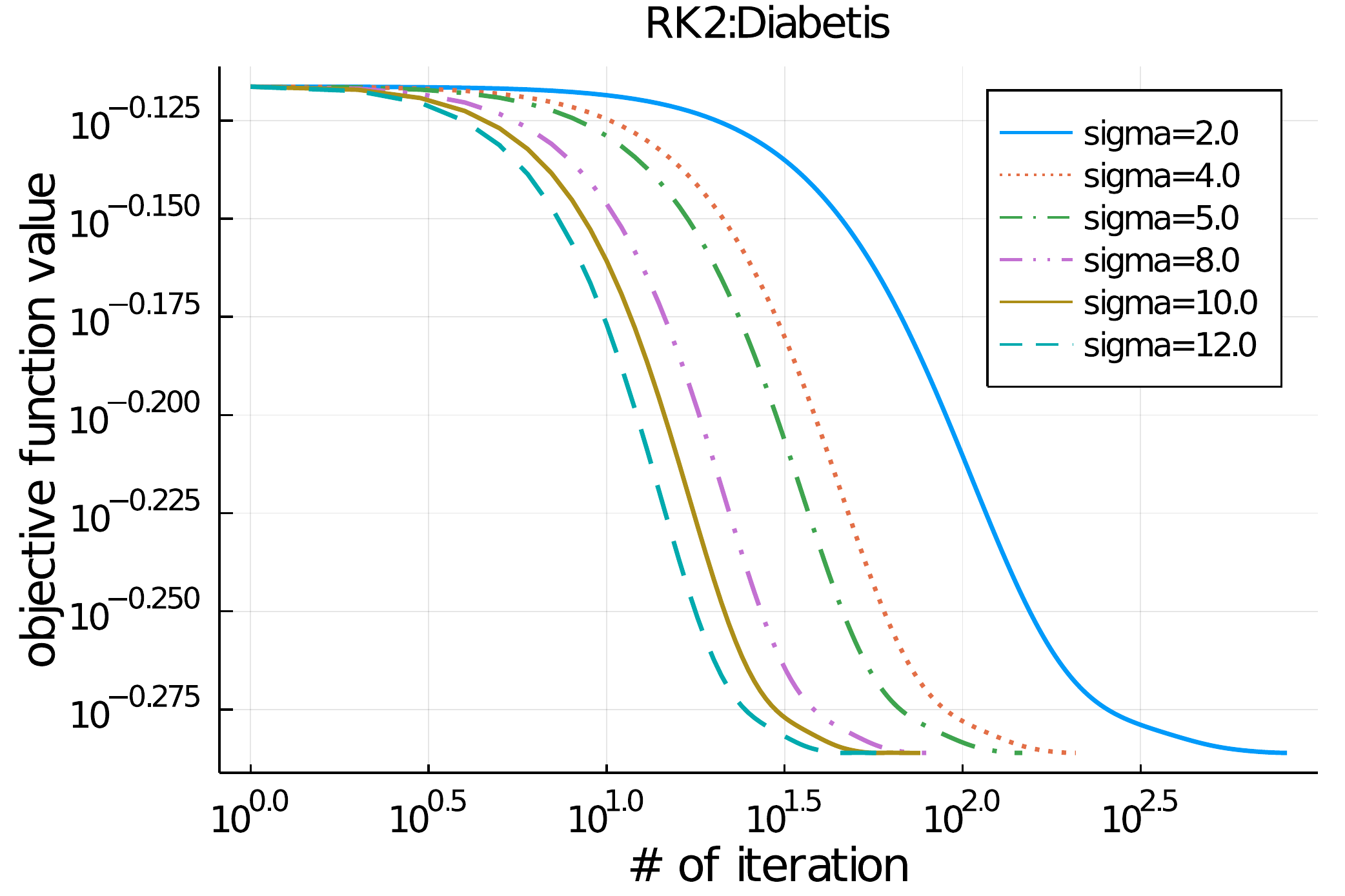}
\includegraphics[scale=0.23]
{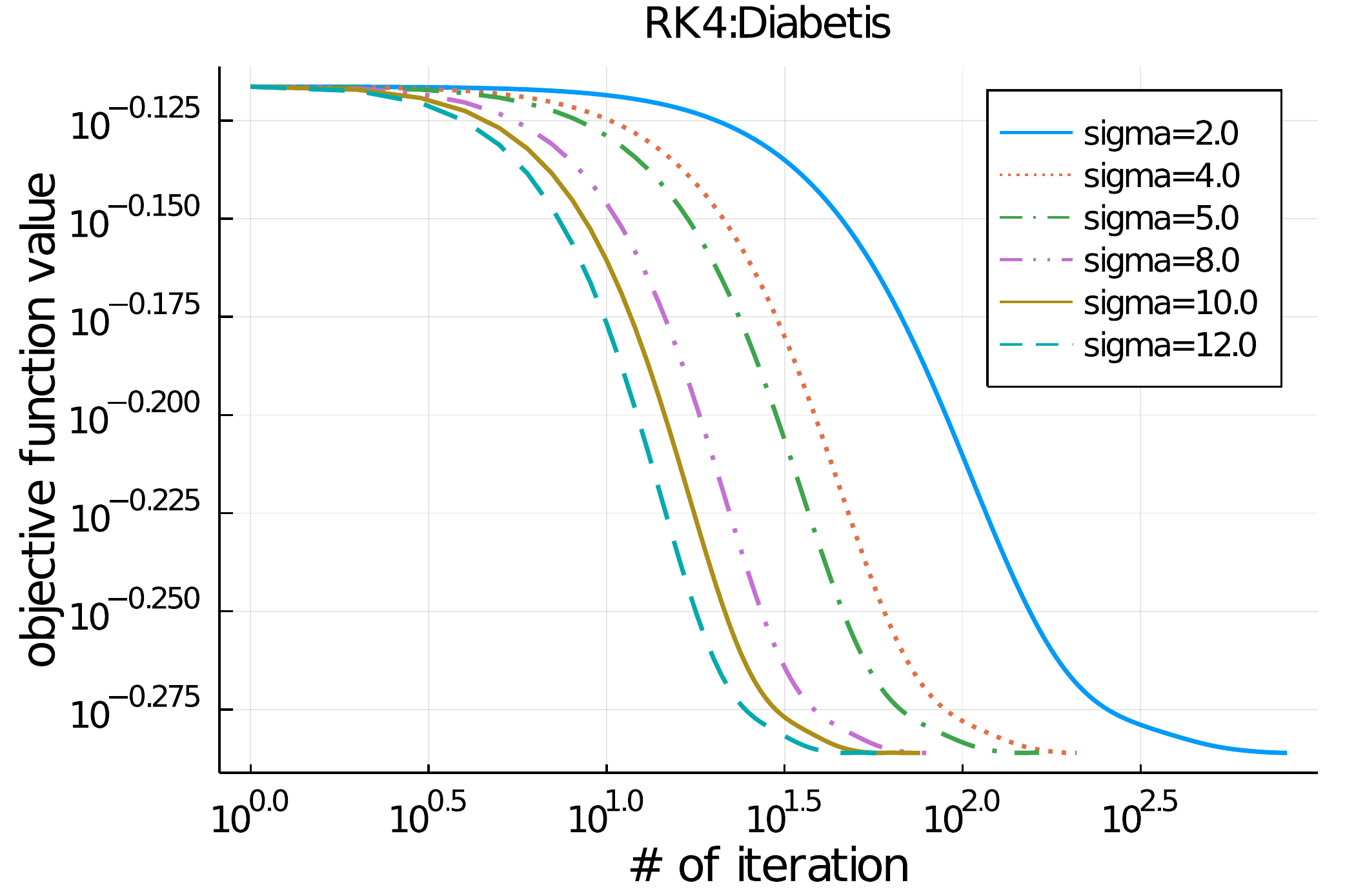}
\includegraphics[scale=0.23]
{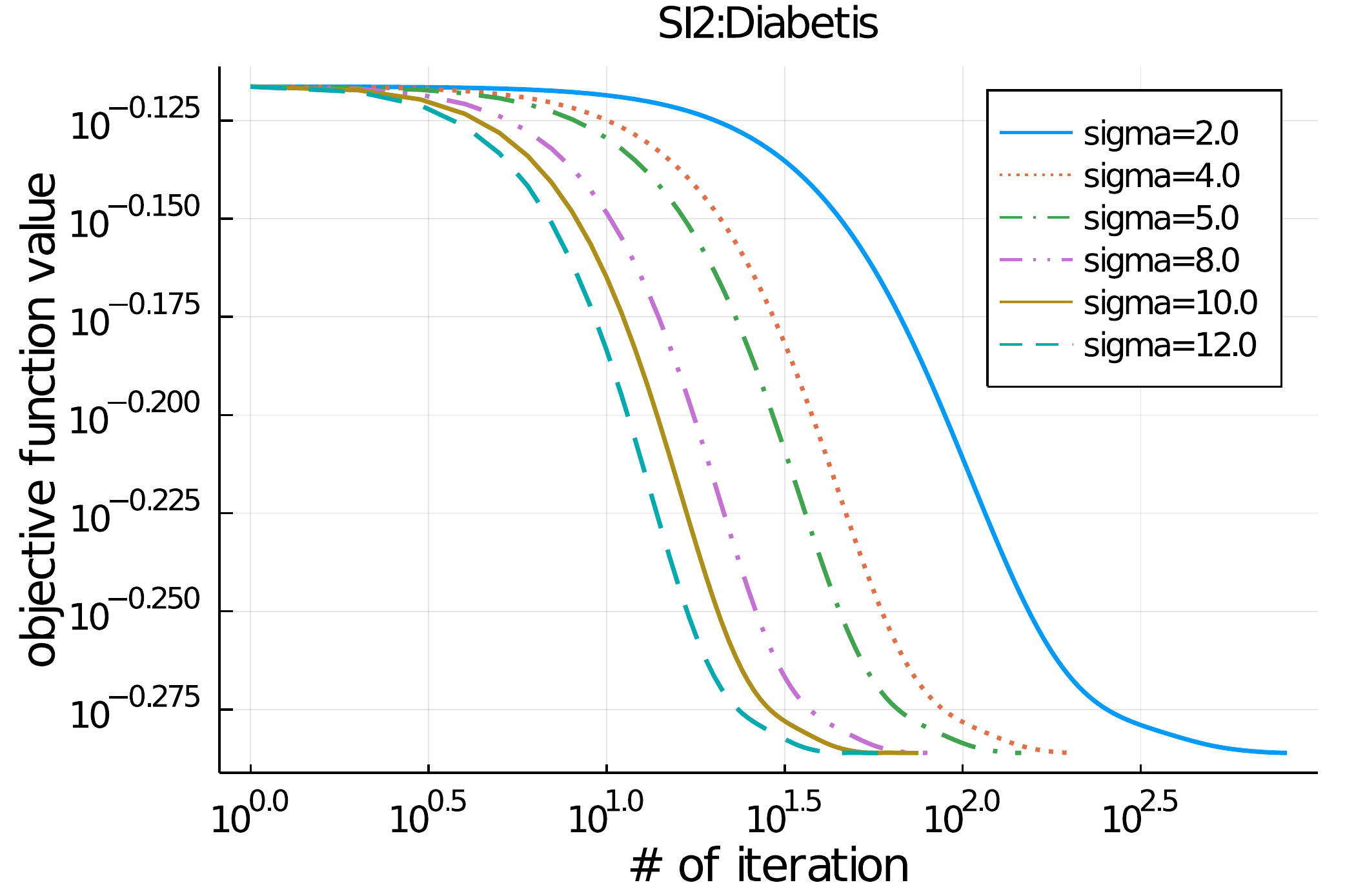}\\
\includegraphics[scale=0.23]
{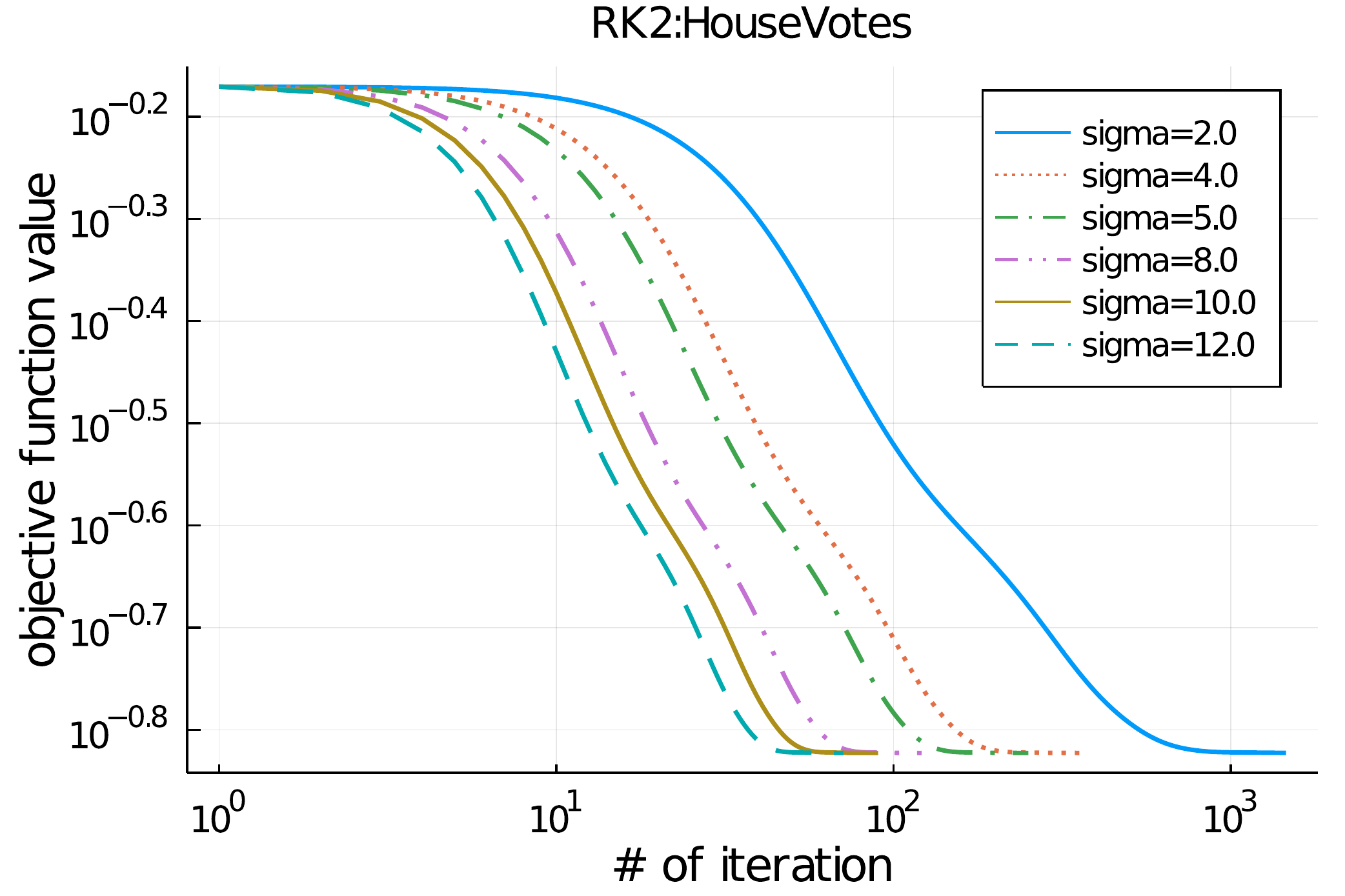}
\includegraphics[scale=0.23]
{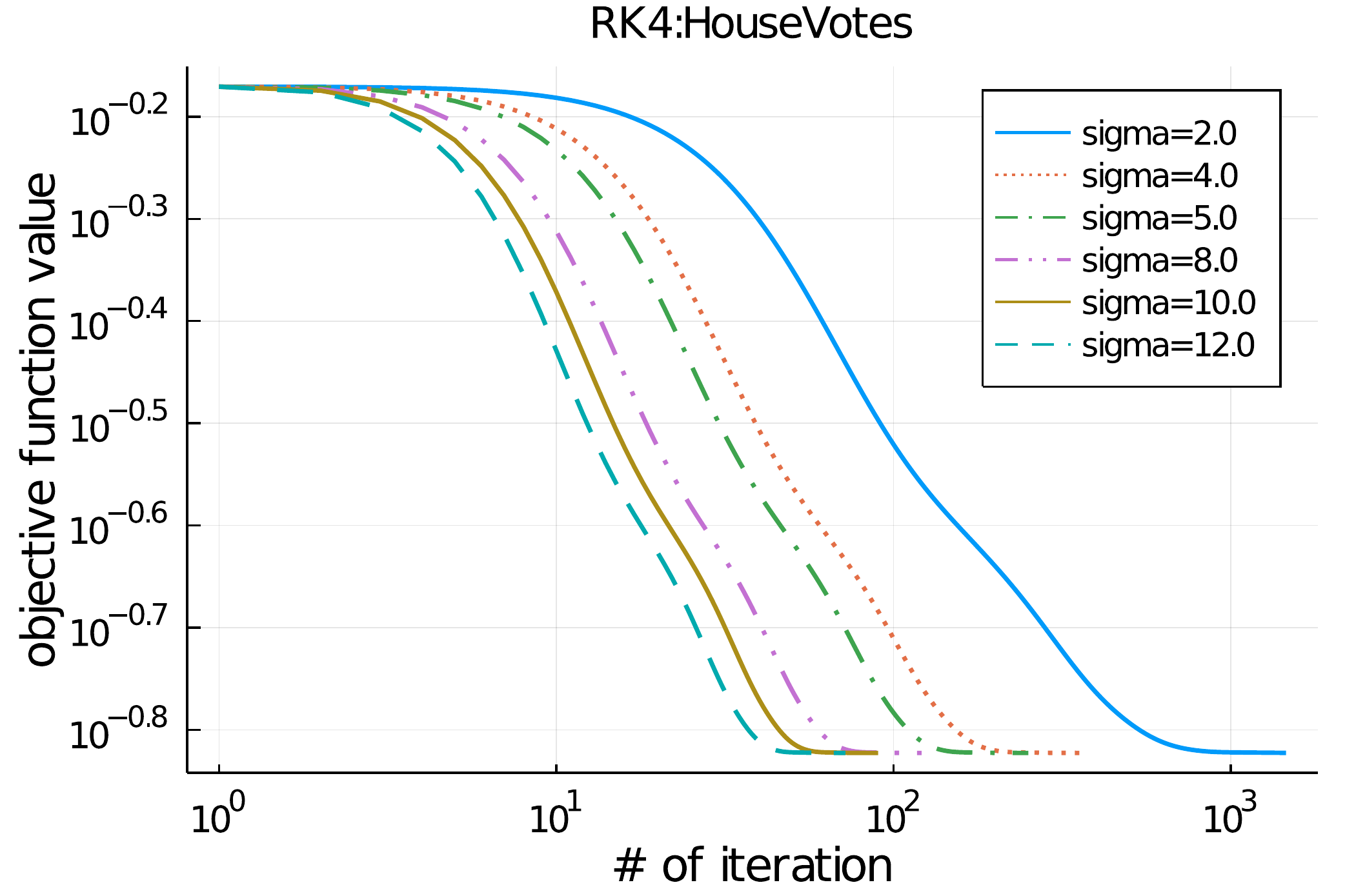}
\includegraphics[scale=0.23]
{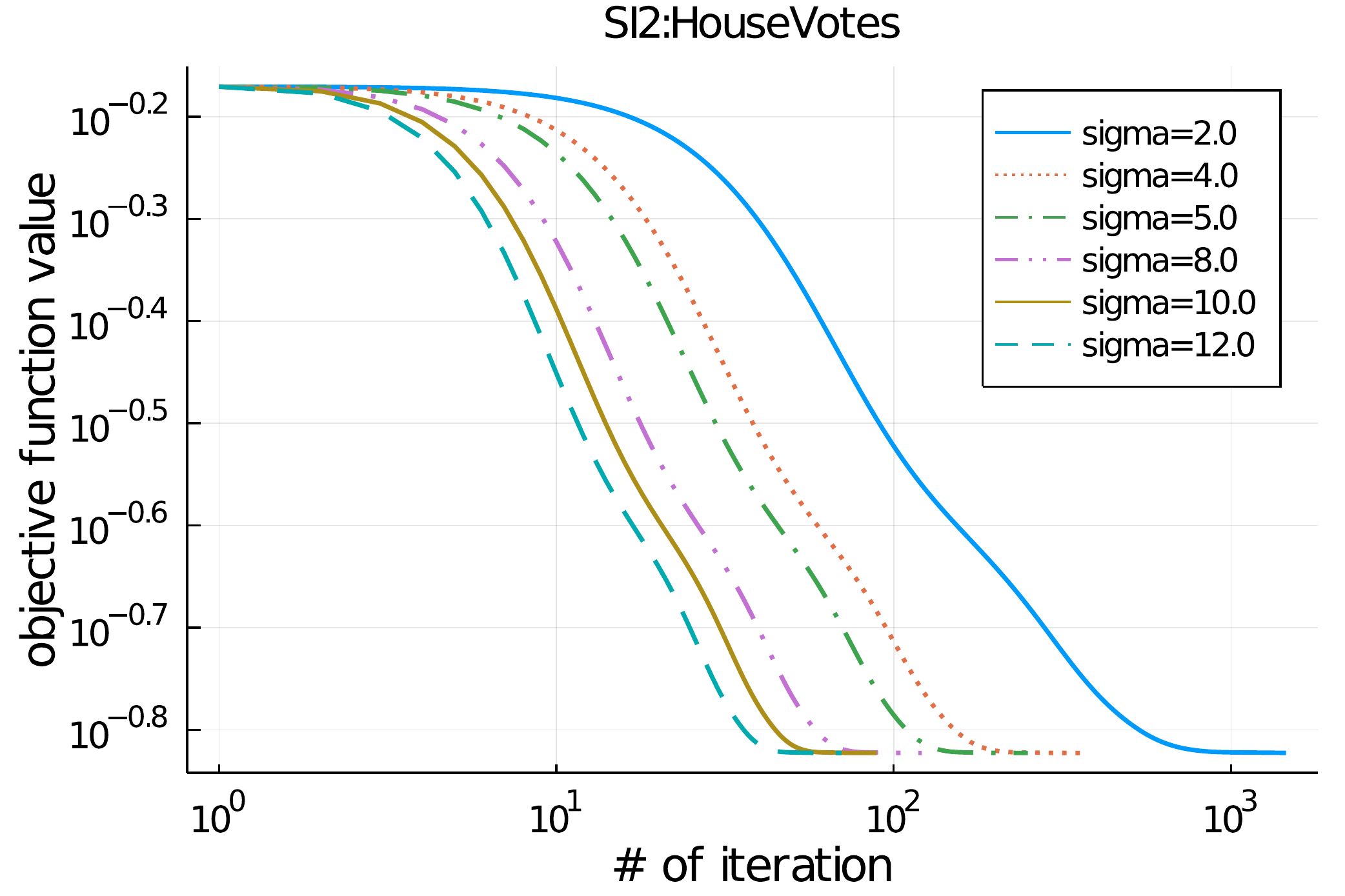}\\
\includegraphics[scale=0.23]
{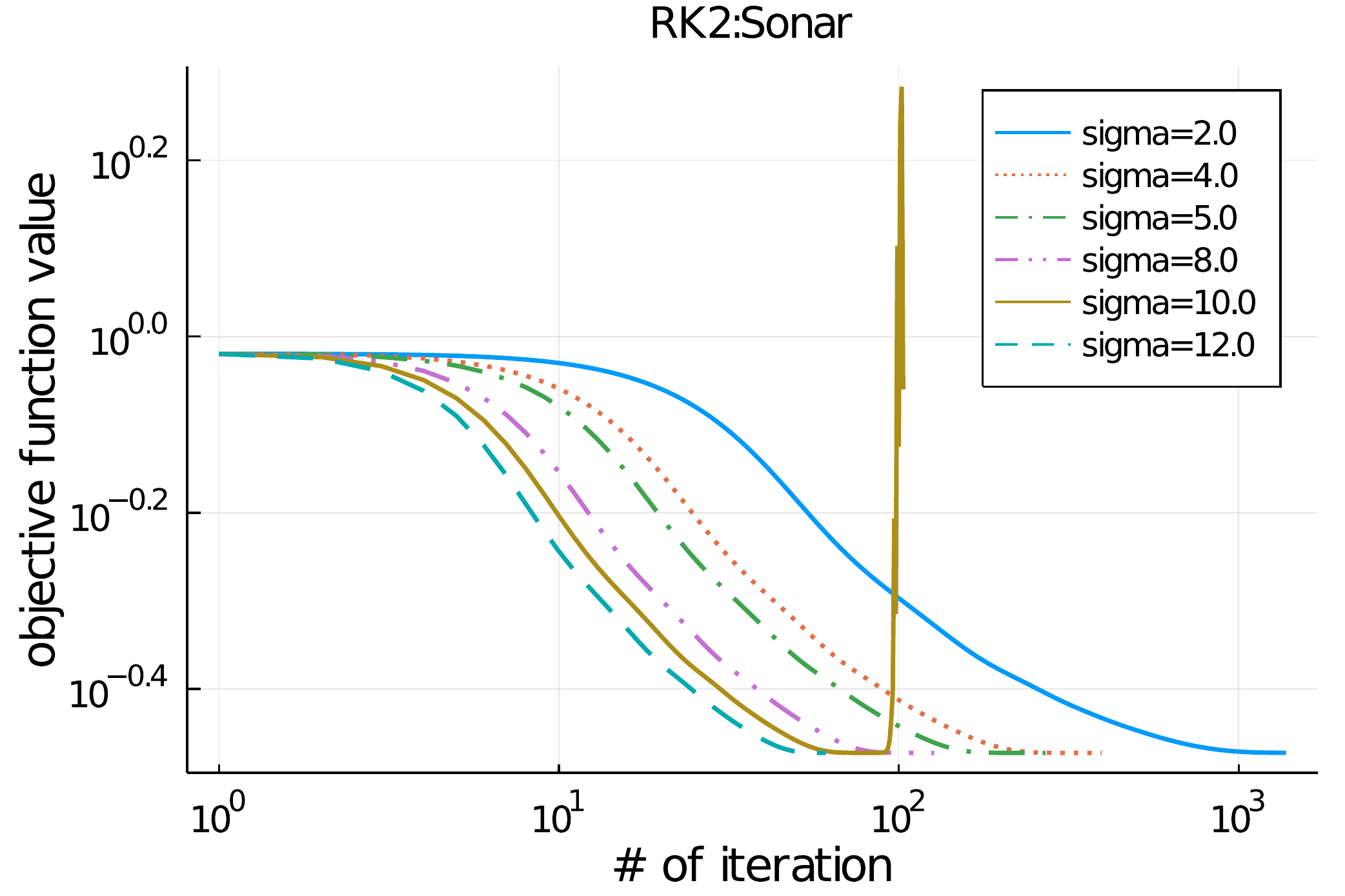}
\includegraphics[scale=0.23]
{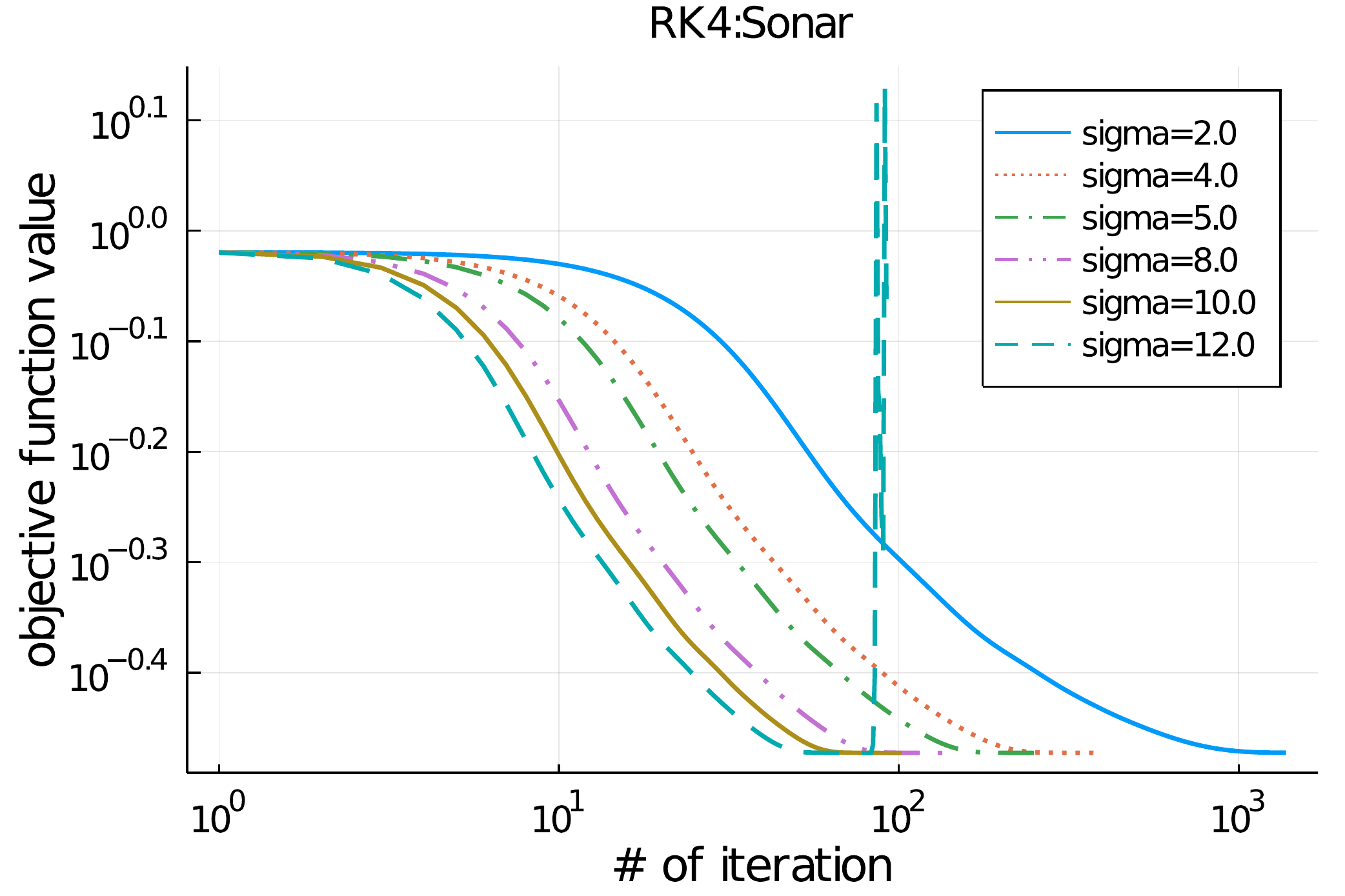}
\includegraphics[scale=0.23]
{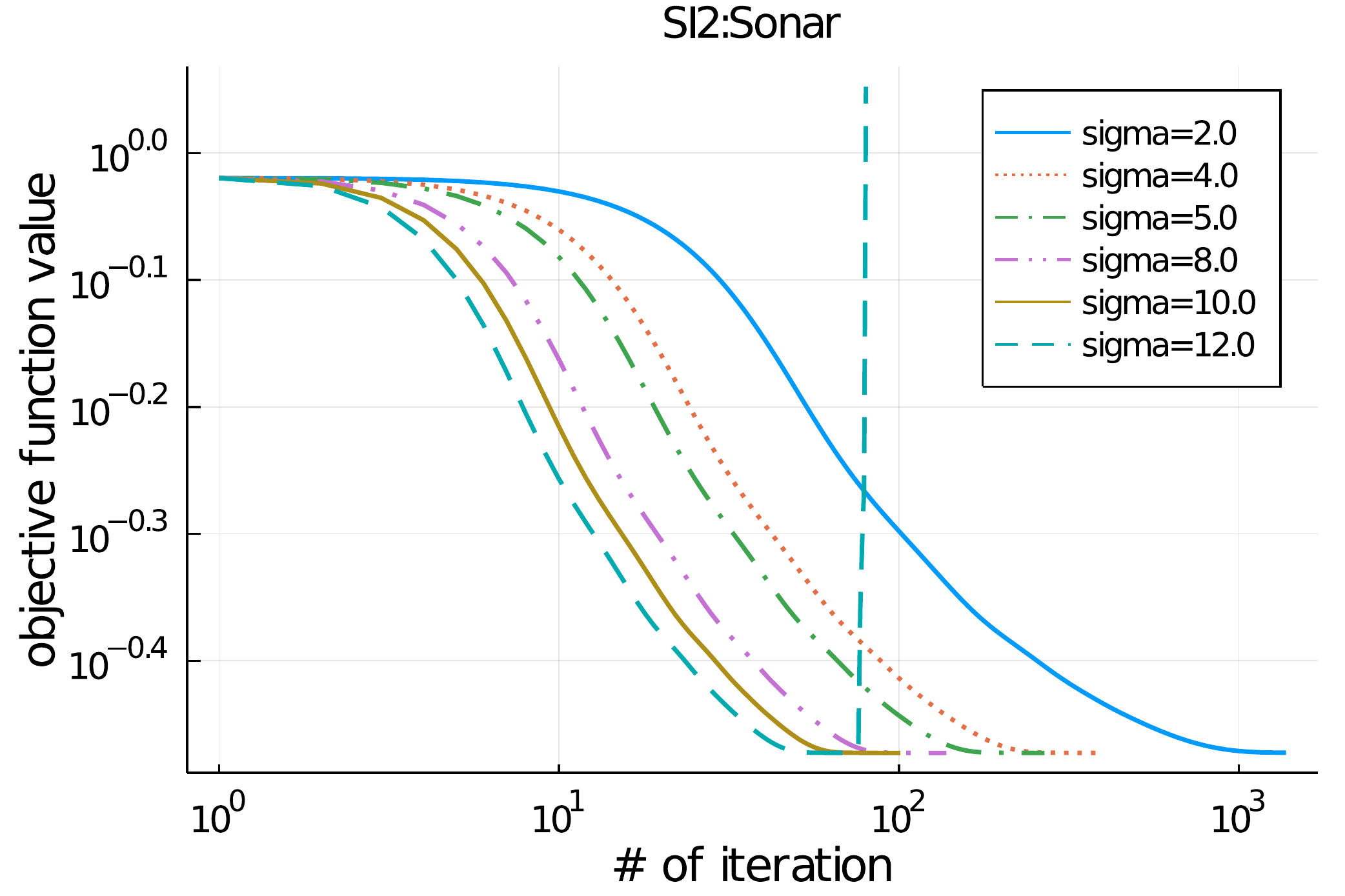}\\
\includegraphics[scale=0.23]
{./figures/RK2_MNIST_l2.pdf}
\includegraphics[scale=0.23]
{./figures/RK4_MNIST_l2.pdf}
\includegraphics[scale=0.23]
{./figures/SI2_MNIST_l2.pdf}
\caption{Objective values along with the iteration of optimization with varying convergence parameter values of $\sigma$.}
\label{fig:supp_convergence}
\end{figure*}

We first consider comparison to 2nd and 4th order Runge-Kutta methods (RK2 and RK4)~\citep{Griffiths2010} with different parameter $\sigma$ in the original ODE defined in Eq.~(6) 
in the main text, which controls the convergence speed. Theoretically, by increasing the value of $\sigma$, we can achieve faster convergence rate. In reality, partly due to discretization errors, too large $\sigma$ may cause numerical instability. We fixed the learning rate $\tau=0.01$ for our SI2, RK2 and RK4, and report the convergence behaviors of objective functions for five datasets in Fig.~\ref{fig:supp_convergence}. 

\begin{figure*}[t!]
\centering
\includegraphics[scale=0.23]
{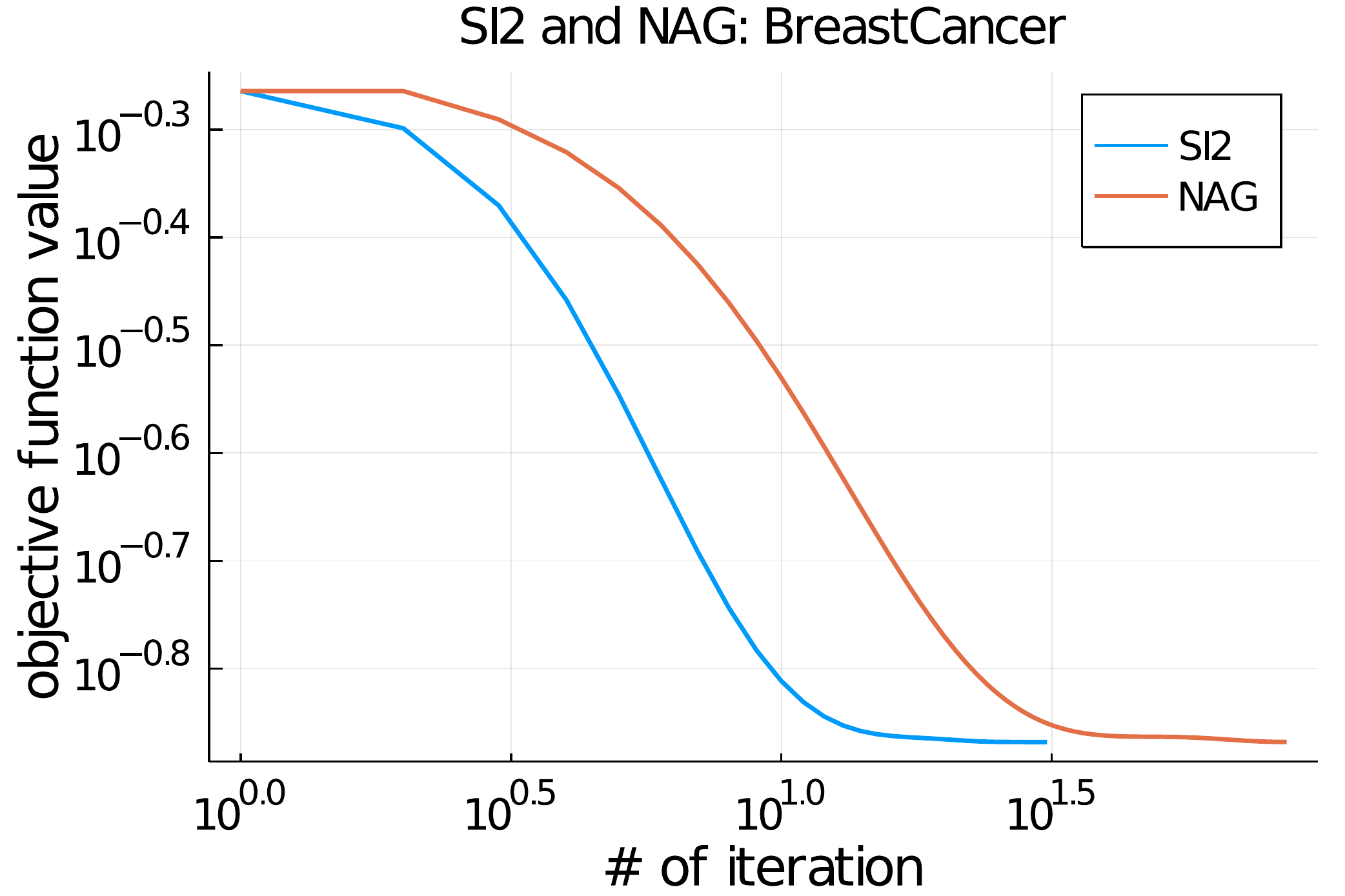}
\includegraphics[scale=0.23]
{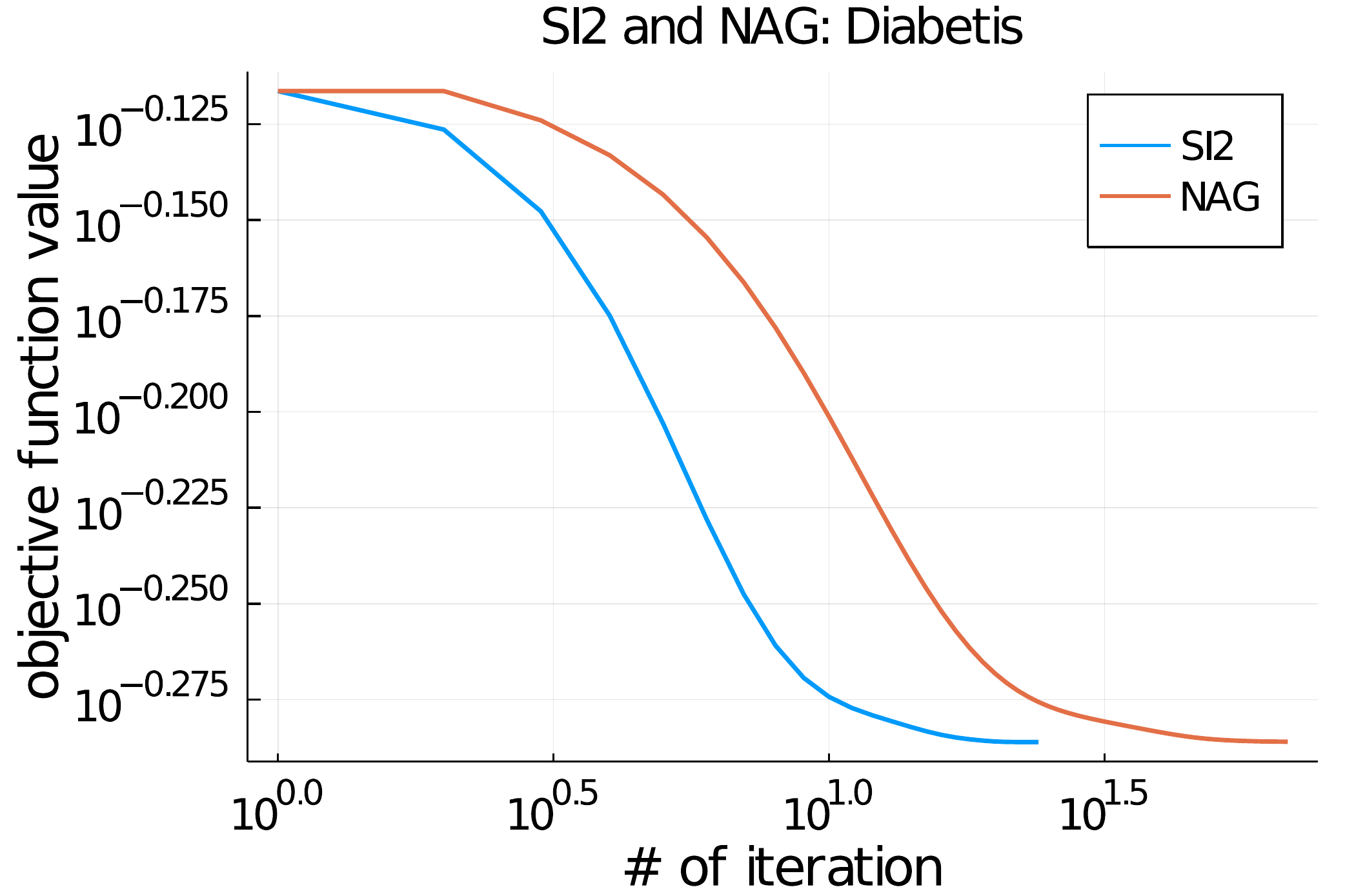}
\includegraphics[scale=0.23]
{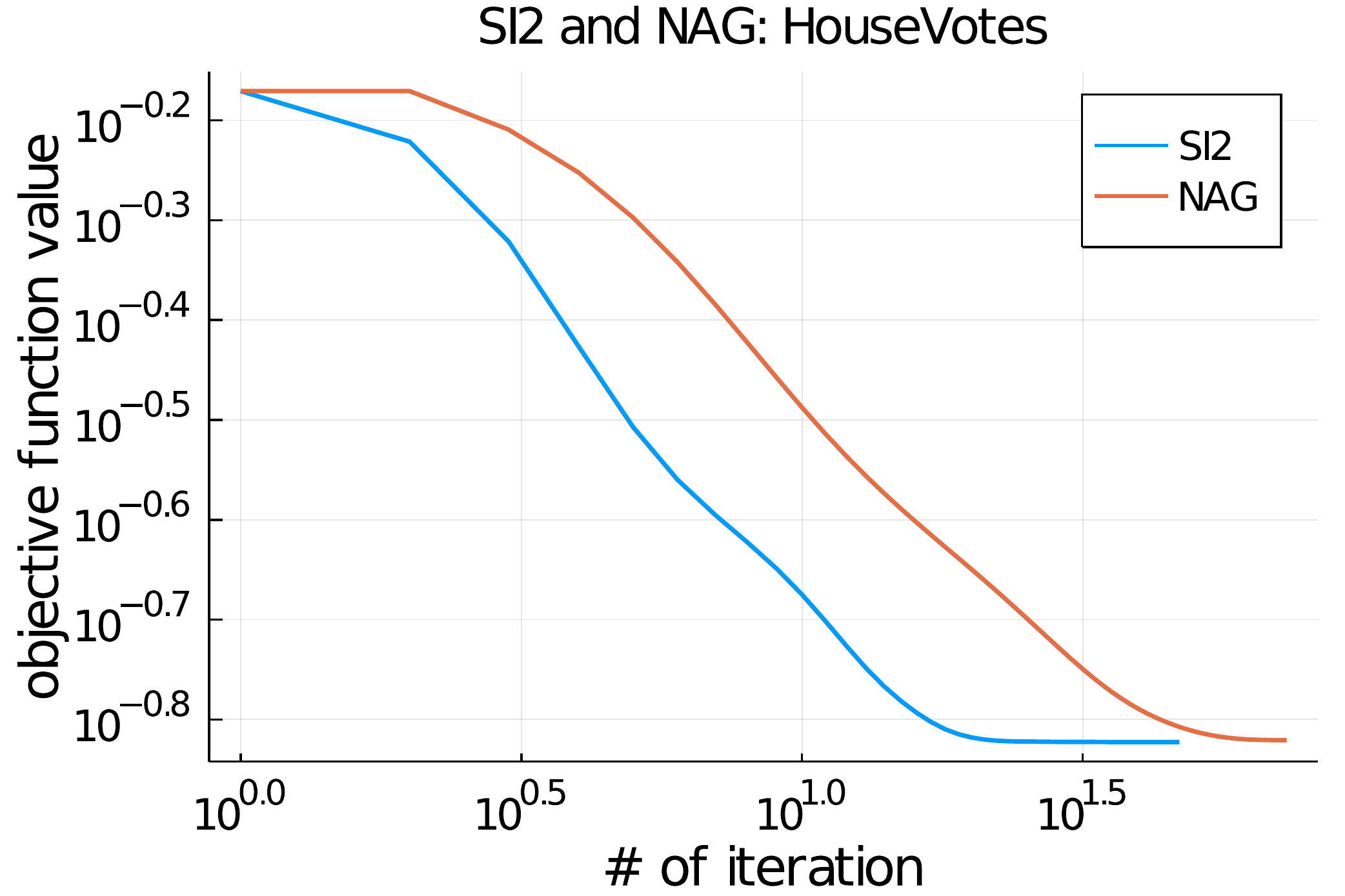}
\includegraphics[scale=0.23]
{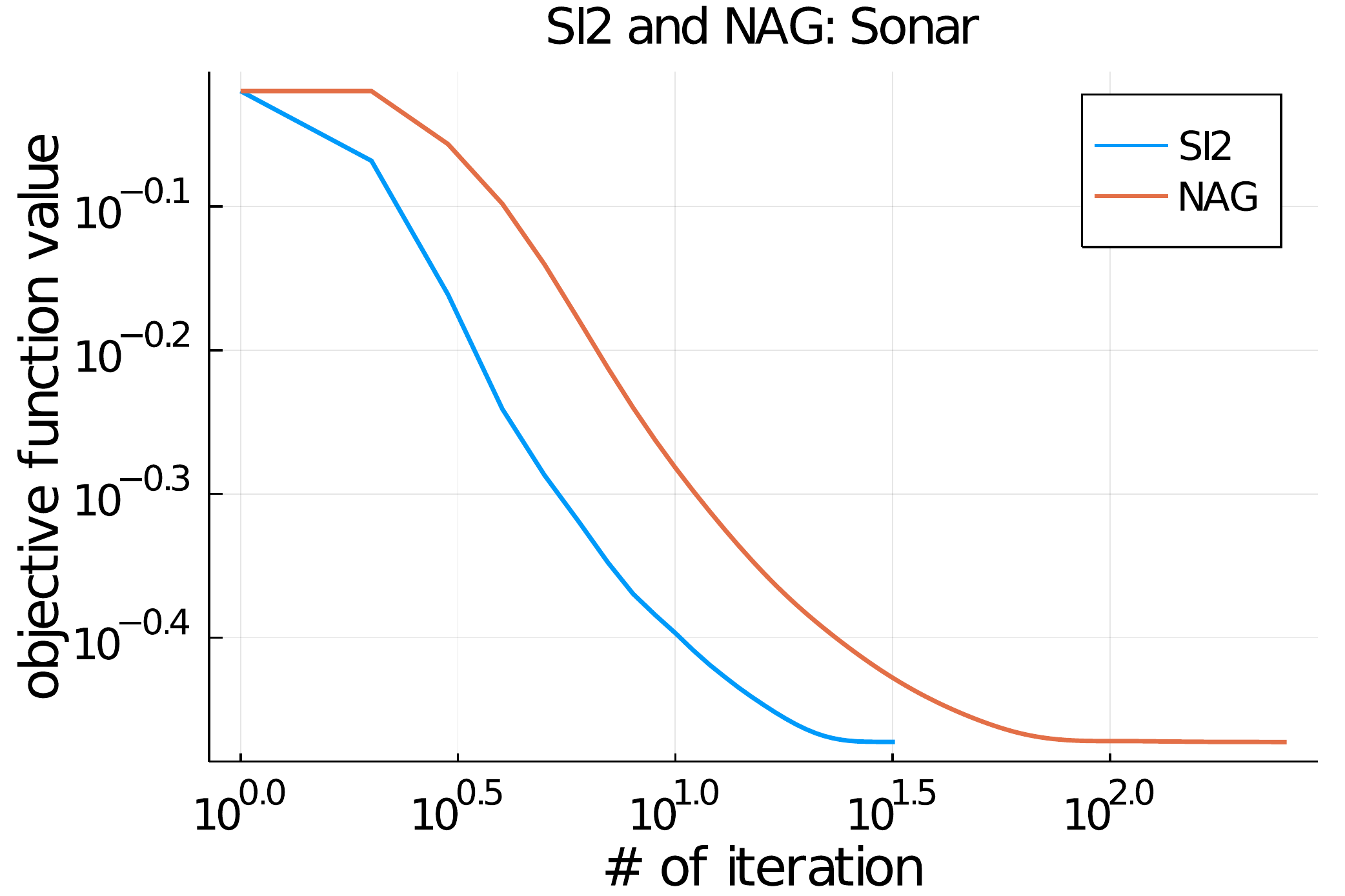}
\includegraphics[scale=0.23]
{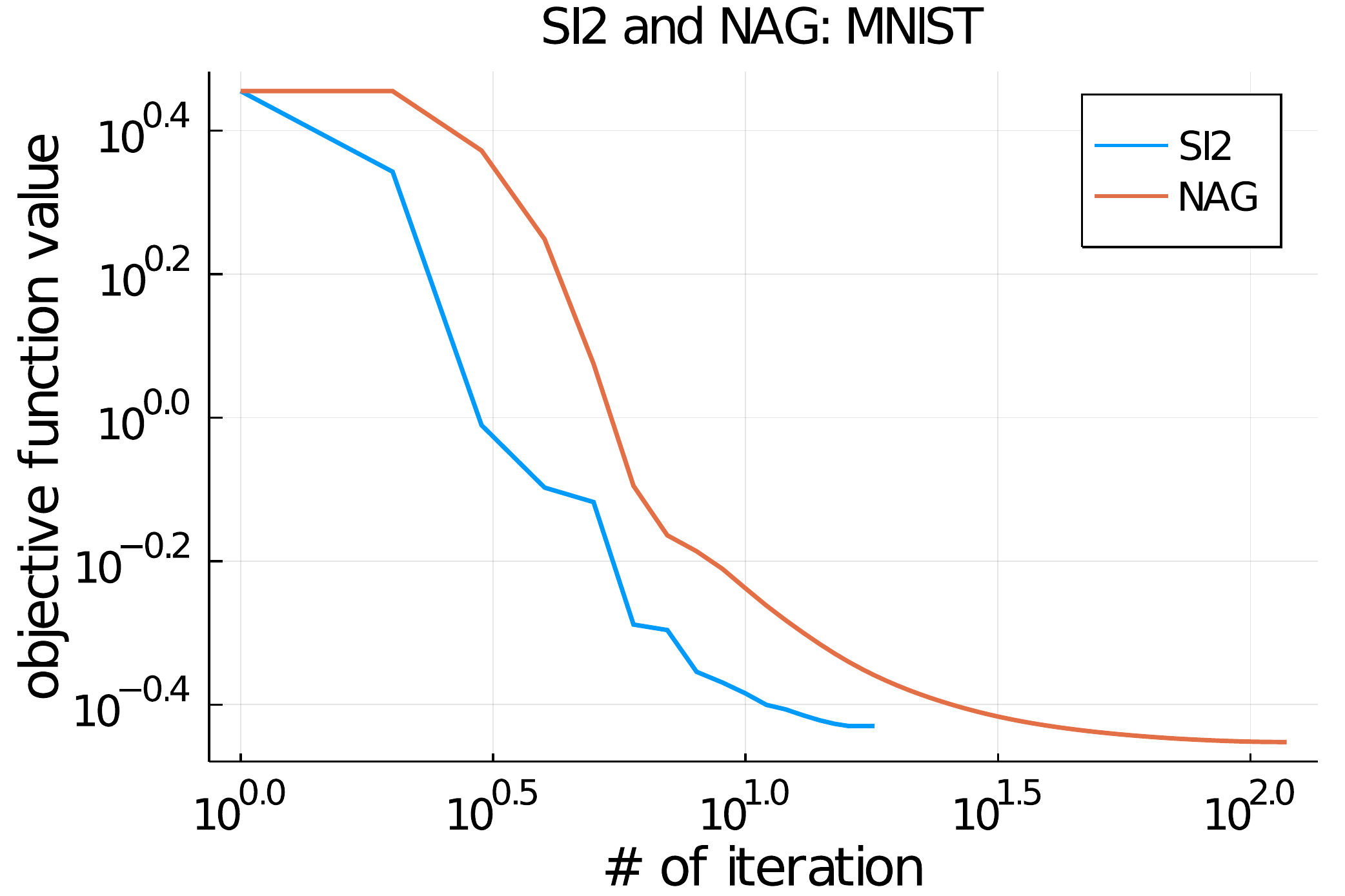}
\caption{Objective values along with the iteration of optimization by SI2 and NAG.}
\label{fig:supp_convergence2}
\end{figure*}

From Fig.~\ref{fig:supp_convergence}, it is seen that in general, with larger $\sigma$, the convergence speed is high. However, too large $\sigma$ cause instability,  particularly for MNIST dataset. SI2 and RK4 are relatively stabler than RK2.
This behavior is partly contributed to the fact that RK2 is based on the lower degree of Taylor expansion and it is deviated from the original ODE compared to RK4.

Now we compare the computational speed of the three methods. Table~\ref{tab:time} shows average of the computational time to reach the stopping criterion (relative difference of the objective function value is less than $10^{-6}$) for RK2, RK4, SI2, and SI2 with backtracking and NAG with backtracking for adjusting step size (see next subsection).
It is seen that SI2 is slightly faster than RK2 or on par for the first four datasets, and faster than RK4. For MNIST, which is the largest size among five datasets, SI2 is significantly faster than other two methods.

We note that SI2, RK2 and RK4 is derived from the same ODE, hence we focus on the convergence behavior and computational cost. 


\subsection{Comparison to NAG method}

We then compare SI2 to NAG described by
\beqa
x^{\,(k)}
&=&y^{\,(k-1)}-s_{\,\N}\,(\nabla f)(y^{\,(k-1)}),
\non\\
y^{\,(k)}
&=&x^{\,(k)}+\frac{k-1}{k+2}\,(x^{\,(k)}-x^{\,(k-1)}), 
\non
\eeqa
where $s_{\,\N}>0$ is a step size parameter, $x^{\,(k)}$ and $y^{\,(k)}$ denote $x\in\mbbR^{\,d}$ and $y\in\mbbR^{\,d}$ at discrete step $k\geq 0$, respectively.

From the results of the previous experiment, we see that for SI2, $\sigma$ less than $8.0$ offers stable results. In this section, the convergence rate parameter $\sigma$ is fixed to $6.0$. For both SI2 and NAG, step size remains to be a tuning parameter. For fair comparison, we adopt the back tracking method for automatically adjust the step size in each iteration. We implemented the momentum restarting mechanism to NAG for stabilizing the performance. 

From Fig.~\ref{fig:supp_convergence2}, it is seen that in general, SI2 requires less iteration for convergence of the objective function value compared to NAG.

From Table~\ref{tab:time}, column SI2(BT) and NAG(BT), it is seen that computational time of SI2 with backtracking is significantly faster than NAG, particularly for MNIST dataset.

\bibliography{reference}
\bibliographystyle{unsrtnat}

\end{document}